\documentclass[11pt]{article}
\input epsf
\usepackage{amsfonts}
\usepackage{amscd}
\usepackage{amsmath,amssymb,amsthm}
\usepackage{amstext}
\usepackage{amscd}
\usepackage{graphicx}
\usepackage[all]{xy}
\usepackage{pstricks,pst-node,pst-tree}

\theoremstyle{plain}
\newtheorem{lem}{Lemma}[section]
\newtheorem*{main}{Theorem 1}
\newtheorem*{smooth_lem}{Lemma \ref{smooth}}
\newtheorem*{cubemoves_thm}{Theorem \ref{cubemoves}}
\newtheorem{theo}[lem]{Theorem}

\newtheorem{coro}[lem]{Corollary}
\newtheorem{prop}[lem]{Proposition}

\theoremstyle{definition}

\newtheorem{definition}[lem]{Definition}

\newtheorem{rem}[lem]{Remark}

\renewcommand{\descriptionlabel}[1]%
       {\hspace{\labelsep}\textsf{#1}}

\begin{document}

\title{Cubulated moves and discrete knots }
\author{ Gabriela Hinojosa, Alberto Verjovsky,
Cynthia Verjovsky Marcotte}

\maketitle
\rightline{\footnotesize Dedicated to Francisco Gonz\'alez Acu\~{n}a,}
\rightline{\footnotesize Fico, on the occasion of his $70^{th}$ birthday.}

\begin{abstract}
In this paper, we prove that given two cubic knots $K_1$, $K_2$ in $\mathbb{R}^3$,  
they are isotopic if and only if one can pass from one to the other by a finite sequence of cubulated moves.  
These moves are analogous to the Reidemeister moves for classical tame knots. We use this fact to
describe a cubic knot in a discrete way, as a cyclic permutation of contiguous vertices of the $\mathbb{Z}^3$-lattice (with some restrictions); moreover, 
we describe a regular diagram of a cubic knot in terms of such cyclic permutations.
\end{abstract}

\section{Introduction}

In \cite{BHV} it was shown that any smooth knot ${\mathbb S}^n\sim{K}^n\subset{\mathbb R}^{n+2}$ can be deformed isotopically into the $n$-skeleton  
of the canonical cubulation of  ${\mathbb R}^{n+2}$ and this isotopic copy is called \emph{cubic knot}. In particular, every classical smooth knot 
$\mathbb S^1\subset{\mathbb R}^3$ is isotopic to a cubic knot.
There are two elementary ``cubulated moves''. The first one (M1) is obtained by dividing each cube of the original cubulation of $\mathbb{R}^3$ into $m^3$ cubes, which means that each 
edge of the knot is subdivided into $m$ equal segments. The second one (M2) consists in exchanging a connected set of edges in a face of the cubulation (or a subdivision 
of the cubulation) with the complementary edges in that face.
If two cubic knots $K_1$ and $K_2$ are such that we can convert $K_1$ into $K_2$ using a finite sequence of cubulated moves then we
say that they are equivalent via cubulated moves, notated $K_1\overset{c}\sim K_2$.\\

This allows us to prove the following:

\begin{main}\label{main}
Given two cubic knots $K_1$ and $K_2$  in $\mathbb{R}^{3}$,
they are isotopic if and only if $K_1$ is equivalent to $K_2$ by a finite sequence of cubulated moves; {\emph{i.e.}}, $K_1\overset{c}\sim K_2$.
\end{main}
Theorem $\ref{main}$ is analogous to the Reidemeister moves of classical tame knots for cubic knots.\\

Since a cubic knot is given by a sequence of edges whose boundaries are 
in the canonical lattice of points with integer coefficients in $\mathbb R^3$,  {\emph i.e.}, the abelian group $\mathbb Z^3$, each knot is determined by a 
cyclic permutation $(a_1,\dots,a_n)$ (with some restrictions), $a_i\in\mathbb Z^3$. In section $5$ we describe a regular diagram of a cubic knot
in terms of such cyclic permutations by projecting onto a plane $P$, such that it is injective when restricted to the $\mathbb{Z}^3$-lattice and the
image of the $\mathbb{Z}^3$-lattice, $\Lambda_P$, is dense. More precisely, the projection of each knot is determined by a cyclic permutation
$(w_1,\dots,w_n)$ (with some restrictions), $w_i\in \Lambda_P$. This fact allows us to develop algorithms to compute some invariants of cubic knots.

\begin{rem} Professor Scott Baldridge has brought to our attention the fact that our main theorem follows also from 
Theorem 1.1 of his paper with Adam Lowrance \cite{BL}. The methods we use are different from theirs 
and we remark that our moves have a natural higher dimensional counterpart. If $K^n\subset{\mathbb R}^{n+2}$ is a cubic knot contained in the $n$-skeleton of the canonical cubulation
${\cal{C}}$ of $\mathbb{R}^{n+2}$, 
then given a cube $Q\in{\cal{C}}$, for a sufficiently fine subdivision of ${\cal{C}}$, the union $A$ of all $n$-dimensional faces of $K$ contained in $\partial{Q}\cong{\mathbb S}^{n+1}$, 
if nonempty, is homeomorphic to a cubulated $n$-disk, hence the closure $A'$ of $\partial{Q} \setminus A$, is also a cubulated disk. The movement (M2) consists in replacing 
$A$ by $A'$. We conjecture that our theorem is also valid in this higher dimensional situation and the proof should be similar to ours.\\

Another motivation for our theorem is to use it to construct dynamically defined wild knots as in  \cite{BHV1}.

\end{rem} 

\begin{rem}
There are several mathematicians working with cubic knots and cube diagrams and their higher-dimensional analogous.  
For example, Ben McCarty was the first to use cube knots to give the easiest known proof that the left-hand and right-hand 
trefoils are not isotopic. Also there are several papers which consider lattices and knots. 
See \cite{mccarty1}, \cite{mccarty2},\cite{mccarty3},\cite{mccarty4},\cite{mccarty5}.

\end{rem}

\section{Cubulations for $\mathbb{R}^{n+2}$}

A {\it cubulation} of $\mathbb{R}^{n+2}$ is a decomposition of $\mathbb{R}^{n+2}$ into a collection $\cal C$ of $(n+2)$-dimensional
hypercubes such that any two of its hypercubes are
either disjoint or meet in one common face of dimension $k\leq n+1$. This provides  $\mathbb{R}^{n+2}$ with the structure of a cubic
complex.\\

In general, the category of cubic complexes and cubic maps is similar to the simplicial category. The main
difference consists in using cubes instead of simplexes. In this context,
a cubulation of a manifold is specified by a cubical complex PL-homeomorphic to the
manifold (see \cite{dolbilin}, \cite{funar}, \cite{matveev}). \\

The {\it canonical cubulation} $\cal C$ of $\mathbb{R}^{n+2}$ is the decomposition of $\mathbb{R}^{n+2}$ into hypercubes which are the images of the unit cube
$I^{n+2}=\{(x_{1},\ldots,x_{n+2})\,|\,0\leq x_{i}\leq 1\}$ by translations by vectors with integer coefficients.\\

\begin{definition} The $n$-skeleton of $\cal C$, denoted by $\cal S$, consists of the union of the $n$-skeletons of the cubes in  $\cal C$,
{\emph {i.e.}}, the union of all cubes of dimension $n$ contained in the faces of the $(n+2)$-cubes in  $\cal C$.
We will call $\cal S$ the {\emph{canonical scaffolding}} of $\mathbb{R}^{n+2}$.
\end{definition}

Any cubulation of $\mathbb{R}^{n+2}$ is obtained from the canonical cubulation by applying  to $\cal{C}$ a conformal transformation
$x\mapsto{\lambda{A(x)}+a}$ where $\lambda\neq0$ is a real number, $a\in \mathbb{R}^{n+2}$ and $A\in{SO(n+2)}$.\\

Consider the homothetic transformation $\frak{h}_{m}:\mathbb{R}^{n+2}\rightarrow\mathbb{R}^{n+2}$ given by
$\frak{h}_m(x)=\frac{1}{m} x$,
where $m>1$ is an integer. The set ${\cal{C}}_{\frac1m}:=\frak{h}_{m}(\cal{C})$ is called a {\it subcubulation} or {\it cubical subdivision} of $\cal C$.\\

In \cite{BHV}, we proved the following theorem which is central to this paper.

\begin{theo} Let $\cal{C}$ be the canonical cubulation of $\mathbb{R}^{n+2}$. Let $K\subset \mathbb{R}^{n+2}$ be a smooth knot of dimension $n$.
There exists a knot
$\hat{K}$ continuously isotopic to $K$, which is contained in the scaffolding $\cal{S}$ of $\cal{C}$
of $\mathbb{R}^{n+2}$. The cubulation of the knot $\hat{K}$ admits a subdivision by simplexes and
with this structure the knot is PL-equivalent to the $n$-sphere with its canonical PL-structure.
\end{theo}

\section{Cubulated isotopy}

Recall that a smooth parametrized knot is a smooth embedding $K: \mathbb{S}^n\to \mathbb{R}^{n+2}$. As is usual, 
we will at times identify the embedding with its image, which we call a geometrical knot.\\

Given two smooth parametrized $n$-dimensional knots $K_1$, $K_2:\mathbb{S}^n\hookrightarrow\mathbb{R}^{n+2}$ we say they are smoothly isotopic 
if there exists a smooth isotopy $H:\mathbb{S}^{n}\times\mathbb{R}\rightarrow\mathbb{R}^{n+2}$ such that
$$
H(x,t)= \left\{ \begin{array}{ll}
K_1(x) & \mbox{if $t\leq-1$},\\
K_2(x) & \mbox{if $t\geq 1$},
\end{array} \right.
$$
and $H(\cdot,t)$ is an embedding of $\mathbb{S}^n$ for $t\in\mathbb{R}$.\\

\begin{definition}
We will say 
$J=\{(H(x,t),t)\in\mathbb{R}^{n+3}\,|\,x\in\mathbb{S}^n,\,\,t\in\mathbb{R}\}$ is the {\emph {isotopy cylinder}} of $K_1$ and $K_2$. Note 
that $J$ is a smooth submanifold of codimension two in $\mathbb{R}^{n+3}$. 
\end{definition}

Let $\cal{C}$ be the canonical cubulation of $\mathbb{R}^{n+3}$.
Let $p:\mathbb{R}^{n+3}\hookrightarrow\mathbb{R}$ be the projection onto the last coordinate. 
\begin{definition}
If $M$ is a connected subset of $\mathbb{R}^{n+3}$ such that $p^{-1}(c)\cap M$ or $p|_M^{-1}(c)$ is connected for all $c\in\mathbb{R}$,
we say that $M$  is {\emph{sliced by connected level sets of $p$}}. 
\end{definition}

\noindent Observe that there is no restriction on the dimension of $M$. \\

\noindent Note that $J$ is sliced by connected sets. The goal of this section is to prove that there exists an isotopic copy $J'$ of $J$ contained
in the $(n+1)$-skeleton of the canonical cubulation $\cal{C}$ of $\mathbb{R}^{n+3}$ which is sliced by connected level sets of $p$.\\

We need some preliminary results:

\begin{prop}\label{cactus}
Let $M$ be a closed subset of $\mathbb{R}^{n+3}$, such that   $p|_M$ is a proper function and $p|_M^{-1}(c)$ is non-empty for all $c\in\mathbb{R}$. Assume $M$ is sliced by connected level sets of $p$. Let $\cal{Q}\subset \cal{C}$ be the union of all cubes of the canonical cubulation intersecting $M$.
 Then $\cal{Q}$ is also sliced by connected level sets of $p$.
\end{prop}

\noindent{\it Proof.} Let $c$ be any number in $\mathbb{R}$. Then $c$ belongs to the closed interval $[n,n+1]$ for some number $n \in \mathbb{Z}$.\\

Let $A$ be the intersection $p^{-1}[n,n+1] \cap M$, so $A$ is compact in $\mathbb{R}^{n+3}$ and is sliced by connected level sets.  Assume $A$ is not connected, so
$A=A_1\cup A_2$ where $A_1$, $A_2$ are non-empty, disjoint compact subsets and both $A_1$ and $A_2$ are also sliced by connected level sets of $p$.
Then $p(A_1)\cup p(A_2)=[n,n+1]$, $p(A_1)\cap p(A_2)=\emptyset$ and $p(A_1)$ and $p(A_2)$ are compact sets, hence $[n,n+1]$ is not connected, which is a contradiction. Therefore $A$ is connected.\\

\noindent{\bf Claim 1:} $p^{-1}[n,n+1]\cap {\cal{Q}}$ is connected.\\ 
Suppose
that $p^{-1}[n,n+1]\cap {\cal{Q}}=C_1\cup C_2$, where $C_1$, $C_2$ are non-empty disjoint closed sets. As $A$ is connected, we have that either $A\subset C_1$
or $A\subset C_2$. We assume that $A\subset C_1$. But $A$ intersects each of the hypercubes belonging to $p^{-1}[n,n+1]\cap \cal{Q}$, since each 
hypercube is connected, it follows that $p^{-1}[n,n+1]\cap {\cal{Q}}\subset C_1$, which is a contradiction. Therefore $p^{-1}[n,n+1]\cap \cal{Q}$ is connected.\\

\noindent{\bf Claim 2:} ${\cal{Q}}$ is sliced by connected level sets of $p$.\\
This is a consequence of the above claim and the fact that the cubes intersecting $p^{-1}(c)$ and  $p^{-1}[n,n+1]$ are the same cubes, since $\mathcal{C}$ is the canonical cubulation. 
$\square$\\

Recall that  $K_1$, $K_2\subset \mathbb{R}^{n+2}$ are two smoothly isotopic knots in $\mathbb{R}^{n+2}$ with $H:\mathbb{S}^{n}\times\mathbb{R}\rightarrow\mathbb{R}^{n+2}$ 
the smooth isotopy between them, which gives us the isotopy cylinder
$J=\{(H(x,t),t)\in\mathbb{R}^{n+3}\,|\,x\in K_1\}$. 

\begin{lem}
The isotopy cylinder $J$ is sliced by connected level sets of $p$.
\end{lem}

\noindent{\it Proof.} Since $p^{-1} (t)=H(\cdot,t)$ is the image of an embedding of $\mathbb{S}^n$, it is connected. Hence the result follows.  $\square$\\

\begin{lem}\label{connected}
 Let $\cal{Q}_J$ be the union of all cubes which intersect $J$. Then  $\cal{Q}_J$ is sliced by connected level sets of $p$.
\end{lem}

\noindent{\it Proof.}  Note that $p|_J$ is clearly proper, therefore this is a consequence of Proposition \ref{cactus}. $\square$\\

\begin{theo}\label{trace}
The isotopy cylinder $J$ can be cubulated. In other words, there exists an isotopic copy $J'$ of $J$ contained
in the $(n+1)$-skeleton of $\cal{Q}_J$, hence $J'$ is contained in the $(n+1)$-skeleton of the canonical cubulation of $\mathbb{R}^{n+3}$. 
Moreover $J'$ can be chosen to be sliced by connected level sets of $p$.
\end{theo}

\noindent{\it Proof.}  Theorem 2.9 in \cite{BHV} implies that $J$ can be cubulated. More precisely in \cite{BHV} it is shown the following:
let $M$ be a closed tubular neighborhood of $J$, then there
 exists a sufficiently small subcubulation ${\cal C}_{\frac{1}{m}}$
 of ${\cal C}$, such that the union $\cal{Q}_{\partial M}$ of all cubes which intersect $\partial M$, is a bicollar neighborhood of $\partial M$. 
We can deform $\partial M$ to any of its boundary components via an adapted flow, which is nonsingular in $M$ and is transverse to $\partial M$.
Let $N\subset \cal{Q}_{\partial M}$ be this isotopic copy. Observe that $N$ contains an
isotopic copy of $J$, say $\tilde{J}$. It was also shown in \cite{BHV} that $\tilde{J}$ can be deformed into the $n$-skeleton of $\cal{Q}_{\partial M}$. This cubic manifold is contained in the subcubulation    ${\cal C}_{\frac{1}{m}}$, hence 
to obtain $J'$ we rescale back the subcubulation to the original one $\cal C$.\\

Finally, notice that $J'$ is contained in the $n$-skeleton of $\cal{Q}_J$, hence the only obstruction for $J'$ to be sliced by connected level sets arises if 
$\cal{Q}_J$ is not sliced by connected 
level sets, but this contradicts Lemma \ref{connected}. $\square$

\section{Cubulated moves}

We say that a knot $K\subset\mathbb{R}^3$ is a {\it cubic knot}, if $K$ is contained in the scaffolding $\cal{S}$ of the 
canonical cubulation $\cal{C}$ of $\mathbb{R}^{3}$.\\

\begin{definition} The following are the allowed {\emph{cubulated moves}}:
\begin{description}
\item[{\bf{M1}}] {\emph{Subdivision:}} Given an integer $m>1$, consider the subcubulation  ${\cal{C}}_m$ of $\cal{C}$. Since ${\cal{C}}\subset {\cal{C}}_{m}$, then
$K$ is contained in the scaffolding ${\cal{S}}_m$ (the 1-skeleton) of ${\cal{C}}_m$ and as a cubic complex, each edge of $K$ is subdivided into $m$ equal
edges.
 
\item[{\bf{M2}}] {\emph{Face Boundary Move:}} Suppose that $K$ is contained in some subcubulation ${\cal{C}}_m$ of  
the canonical cubulation ${\cal{C}}$ of $\mathbb{R}^{3}$. 
Let $Q\in {\cal{C}}_m$ be a cube such that
$A=K\cap Q$ contains an edge. We can assume, up to applying the elementary move $M1$, that $A$ consists of either one
or two edges that are connected and are part of the boundary of a face
$F\subset Q$. Thus $A$ is an arc contained in the boundary of
$F$ and $\partial F$ is divided by $\partial A$ into two cubulated arcs. One of them is $A$ and we denote the other by $A'$. Observe
that both arcs share a common boundary. The move  consists in replacing $A$ by $A'$. 
\end{description}
\end{definition}

\begin{rem}
 Notice that the face boundary move can be extended to an ambient isotopy of $\mathbb{R}^3$.
\end{rem}

\begin{definition}
Given two cubic knots $K_1$ and $K_2$ in $\mathbb{R}^3$. We say that
$K_1$ is {\emph{equivalent}} to $K_2$ {\emph{by cubulated moves}}, denoted by $K_1\overset{c}\sim K_2$, if we can transform $K_1$ to $K_2$ by a finite number
of cubulated moves. 
\end{definition}

\subsection{Main theorem}

We are now ready to prove Theorem \ref{main}.

\begin{main}
Given two cubic knots $K_1$ and $K_2$  in $\mathbb{R}^{3}$,
 then they are isotopic if and only if $K_1$ is equivalent to $K_2$ by cubulated moves; {\emph{i.e.}}, $K_1\overset{c}\sim  K_2$.
\end{main}

\noindent{\it Proof.}  
First, note that if $K_1$ and $K_2$ are equivalent by cubulated moves, then these knots must clearly be isotopic. Hence, what 
remains to be proved is that two cubic knots that are isotopic must also be equivalent by cubulated moves.\\

Our strategy is as follows. First, for $i\in \{1,2\}$, we will smooth each $K_i$ to obtain $\widetilde{K}_i$, and then 
cubulate these two knots to obtain $K'_i$ in such a way that a) $K_i\overset{c}\sim K'_i$ and b) we can show 
$K'_1\overset{c}\sim  K'_2$ are equivalent by cubulated moves.\\

Given a cubic knot $K$, there exists a smooth knot $\widetilde{K}$ isotopic to $K$ such that $\widetilde{K}$ is 
${\cal{C}}^0$-arbitrarily close to $K$. This is because
we can round the corners at the vertices of $K$ in an arbitrarily small neighborhoods of them (see \cite{douady}). \\

Let $J$ be the isotopy cylinder (defined above) of $\widetilde{K}_1$ and $\widetilde{K}_2$.
Then $J$ is a smooth submanifold of codimension two in $\mathbb{R}^{4}$. 
By Theorem \ref{trace}, there exists an isotopic copy of $J$, say $J'$, contained in the $2$-skeleton of the canonical 
cubulation $\cal{C}$ of $\mathbb{R}^{4}$. Recall that $J'$ is sliced by connected level sets of $p$. 
Furthermore, note that there exist integer numbers 
$m_1$ and $m_2$ such that $p^{-1}(t)\cap J' = K'_1$ for all $t\leq m_1$ and
$p^{-1}(t)\cap J' = K'_2$ for all $t\geq m_2$, where $K'_1$ and $K'_2$ are cubic knots which are isotopic to $\widetilde{K}_{1}$ and $\widetilde{K}_{2}$, respectively.\\

Now, we will use the following results which will be proved in Sections 4.1.1 and 4.1.2, respectively.\\

\begin{lem}\label{smooth} 
Given a cubic knot $K$ we can choose a small cubulation ${\cal{C}}_{\frac{1}{m}}$  fine enough that 
$N(K)=\cup\{Q\in {\cal{C}}_{\frac{1}{m}}\,|\,Q\cap K\neq\emptyset\}$ is a closed tubular neighborhood of $K$ and
$Q\cap K$ is equal to either a vertex, one edge, or two edges sharing a vertex (neighboring edges). We can also choose $\widetilde{K}$ isotopic to $K$ 
such that $\widetilde{K}$ is 
${\cal{C}}^0$-arbitrarily close to $K$ and $\widetilde{K}\subset \mbox{Int}(N(K))$. Let 
$K'$ be an isotopic copy of $K$ contained in $\partial N(K)$, then $K\overset{c}\sim  K'$; {\emph{i.e.}}, we can
go from $K$ to $K'$ by a finite sequence of cubulated moves.
\end{lem}

\begin{rem}\label{conditions}
\begin{enumerate}
\item The existence of $K'$ is proved on Theorem 3.1 in \cite{BHV} (see also the proof of Theorem \ref{trace}).
\item We may assume, using a subdivision move if necessary, that the intersection $Q\cap K$ is equal to either a vertex, one edge, or two neighboring edges. 
\end{enumerate}
\end{rem}

\begin{theo}\label{cubemoves} 
Given two cubic knots $K_1$ and $K_2$, we obtain $K'_{1}$ and $K'_{2}$ as in Lemma \ref{smooth}. Then
there exists a finite sequence
of cubulated moves that carries $K'_{1}$ into  $K'_{2}$. In other words, $K'_{1}$ is equivalent to  $K'_{2}$ by cubulated moves.
\end{theo}

Thus by Lemma \ref{smooth}, there exists a finite sequence of cubulated moves that carries $K_{1}$ into $K'_{1}$ and also a finite
sequence of cubulated moves that carries $K_{2}$ into $K'_{2}$. By Theorem \ref{cubemoves},
 there exists a finite sequence of cubulated moves that carries $K'_1$ into $K'_2$. As a consequence there exists a finite sequence
of cubulated moves that converts $K_{1}$ into ${K}_{2}$. $\square$
  
\subsubsection{Cubic case}

Let $K$ be a cubic knot. We can choose a cubulation ${\cal{C}}_{\frac{1}{m}}$  fine enough such that the subcollection
${\cal{N}}_K=\{Q\in {\cal{C}}_{\frac{1}{m}}\,|\,Q\cap K\neq\emptyset\}$ satisfies that
$N(K)=\cup_{Q\in {\cal{N}}_K} Q$ is a closed tubular neighborhood of $K$ and
$Q\cap K$ is equal to either a vertex, one edge, or two neighboring edges (see algo Remark \ref{conditions}). We can also choose $\widetilde{K}$ isotopic to $K$ 
such that $\widetilde{K}$ is 
${\cal{C}}^0$-arbitrarily close to $K$ such that for any cube $Q\in {\cal{N}}_K$ the intersection $\widetilde{K}\cap Q$ is nonempty, 
and $\widetilde{K}\subset \mbox{Int}(N(K))$. 
Observe that $N(K)$ is also a closed tubular neighborhood of $\widetilde{K}$. By Theorem 3.1 in \cite{BHV},
there exists an isotopic copy $K'$ of $\widetilde{K}$ contained in the
1-skeleton  of  $\partial N(K)$. \\

In this section, we will prove that there exists a finite sequence of cubulated moves that carries $K$ into $K'$. 

\begin{smooth_lem}
There exists a finite sequence of cubulated moves that carries $K$ into $K'$.
\end{smooth_lem}

\noindent{\emph{Proof.} 
Observe the following. The knots $K$ and $K'$ are isotopic, both are contained in the 1-skeleton of $N(K)$, but $K\subset \mbox{Int}(N(K))$  and
$K'\subset\partial N(K)$. Our goal is to construct a surface ${\cal {F}}$  such that it will be contained in the 2-skeleton of $N_k$ and its boundary 
will consist of two connected components, namely $K$ and $K'$.\\

Let $B=\{Q\subset{\cal{N}}_k \,|\,Q\cap K\neq\emptyset,\,\,Q\cap K'\neq\emptyset\}$. \\

Since all cubes in ${\cal{N}}_K$ intersect $K$, $B$ consists
of a finite number of cubes, say $m$, such that all of them intersect $K'$. By orienting $K'$ we can enumerate the cubes in $B$ in such a way that
consecutive numbers belong to neighboring cubes (cubes
sharing a common face) and the $m^{th}$-cube is next to the first one.
To construct the surface ${\cal {F}}$, we will look at all possible cases of $Q_j\in B$ 
and find pieces $F_j$, which will be consist of either one face or two neighboring faces of $Q_j$ and whose union will be ${\cal {F}}$. 
The boundary of these $F_j$'s will 
intersect both $K$ and $K'$, hence faces corresponding to neighboring cubes will share a common edge.  \\

\noindent{\bf Claim 1:} Let $Q_j\in B$. If $Q_j\cap K$ consists of a vertex and $Q_j\cap K'$ consists of either a vertex or an edge, then 
$K\cup K'\subset \cup_{Q_i\in B\setminus\{Q_j\}} Q_i$. In other words, $Q_j$ is superfluos. \\
\noindent {\it Proof of Claim 1.} We will look at all cases.
\begin{description}
\item [$(a)$] Each intersection ($Q_j\cap K$ and $Q_j\cap K'$) consists of a vertex. This implies that each knot passes through two edges sharing the
corresponding vertex and these edges are contained in cubes belonging to $B\setminus\{Q_j\}$. 
\item [$(b)$] $Q_j\cap K$ consists of a vertex $v$ and $Q_j\cap K'$ consists of an edge $e$ or vice versa (see Figure \ref{I1}). 
Suppose that both intersections are contained in the same face $F\subset Q_j$ (see Figure \ref{I1}(a)).
Notice that $K$ passes through two edges sharing the vertex $v$. These two edges belong to neighbor cubes of $Q_j$. By construction of both $N(K)$ and $B$, we 
can assume that these cubes belong to $B$.  This implies that
one of these two neighbor cubes contains both $v$ and $e$. \\
Now, suppose that the vertex $v$ is opposed to $e$, \emph{i.e.}, there is no edge  of $Q_j$ joining
$v$ to any of the end-points of $e$ (see Figure \ref{I1}(b)). As in the previous case, $K$ passes through two edges sharing the vertex $v$ and 
these two edges belong to neighbor cubes of 
$Q_j$ such that $K'$ intersects them. We have that one of these cubes must contain both $v$ and one of the end-points of $e$. 
However this configuration implies that  $K'$ does not intersect any of the remaining neighbor cubes $Q_r$ of $Q_j$ such that $v\in Q_r$. This is a contradiction, 
hence this configuration
is not possible. Therefore, the claim 1 follows. $\square$
\begin{figure}[h] 
 \begin{center}
 \includegraphics[height=3cm]{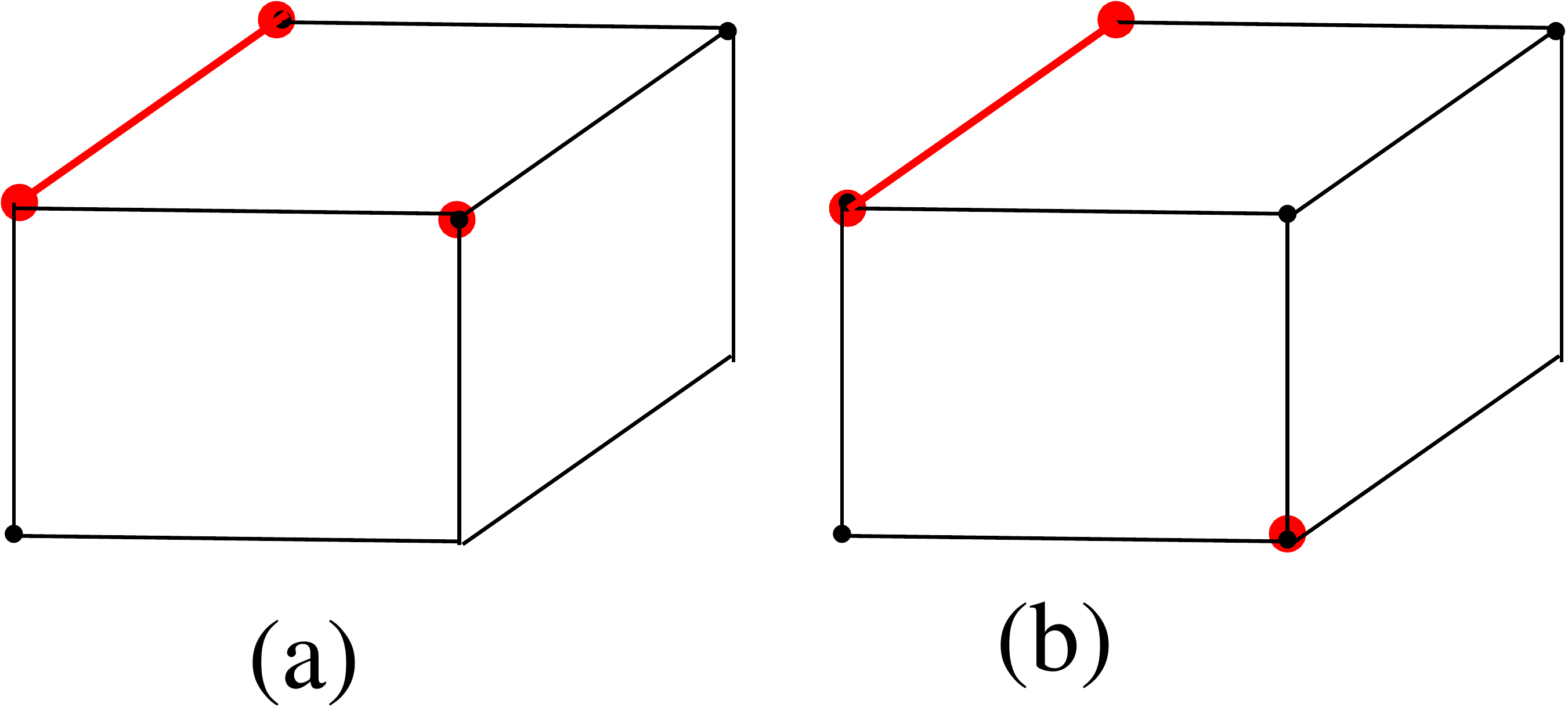}
\end{center}
\caption{\sl $Q_j\cap K$ consists of a vertex $v$ and $Q\cap K'$ consists of an edge $e$.} 
\label{I1}
\end{figure} 

\end{description}

\noindent For simplicity, we will assume that $B$ does not contained superfluous cubes.\\

\noindent {\bf Claim 2:} Let $Q_j\in B$. Then, using face boundary moves if necessary, we can assume that 
$K'\cap Q_j$ is a path consisting of the union of at most three edges.\\
\noindent {\it Proof of Claim 2.} We will prove that $K'\cap Q_j$ can be reduced using face boundary moves if necessary, to a path consisting 
of the union of at most three edges.  
We will consider all possible cases such that $K'\cap Q_j$ consists of the union of at least four edges. Remember that 
$K\cap Q_j$ is equal to either a vertex, one edge or two neighboring edges.

 \begin{description}
 \item [Case 1.] Suppose that $K\cap Q_j$ consists of a vertex $v$. Since $K\cap K'=\emptyset$, it follows that $K'\cap Q_j$ 
can be contained in either one, two or three faces of $Q_j$ and
these faces do not contained $v$. By a combinatorial analysis, we have that $K'\cap Q_j$ is a path consisting of at most 6 edges. 
\begin{description}
\item [$(a)$] $K'\cap Q_j$ consists of a six edges path. In this case,  three edges are contained in a face $F_1$, two disjoint edges lie on a neighbor face $F_2$ of
$F_1$ and one
edge lies on a neighbor face $F_3$ of $F_1$(see Figure \ref{I2}). Now, we apply the face boundary move (M2) on $F_1$ to obtain a four edges path, where three of them
are contained in $F_2$. Next, 
we apply again an (M2)-move on $F_2$ to obtain that  $K'\cap Q_j$ can be 
reduced to a two edges path.
\begin{figure}[h] 
 \begin{center}
 \includegraphics[height=3cm]{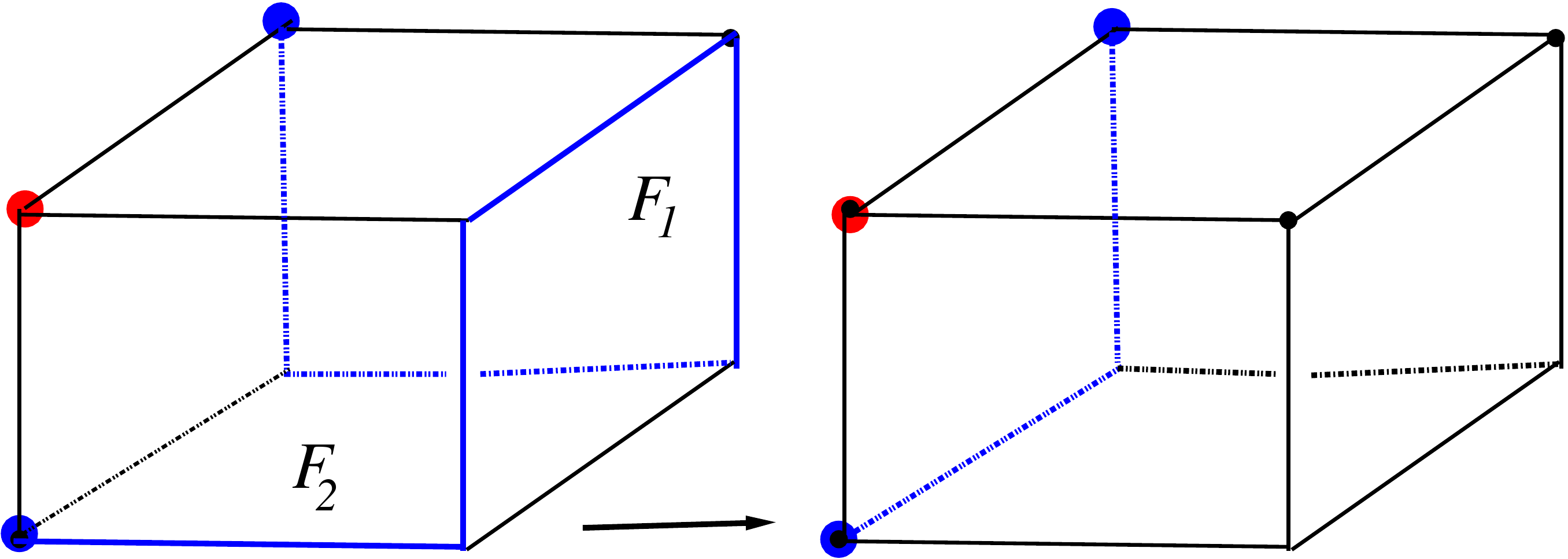}
\end{center}
\caption{\sl $K'\cap Q_j$ consists of a six edges path.} 
\label{I2}
\end{figure} 
\item [$(b)$] $K'\cap Q$ consists of a five edges path. In this case, three edges are contained in a face $F_1$ and the remaining two edges lie on a neighbor face 
$F_2$. Then, applying an (M2)-move on  $F_1$ we obtain that $K'\cap Q_j$ can be reduced to a three edges path (see Figure \ref{I3}).
\begin{figure}[h] 
 \begin{center}
 \includegraphics[height=3cm]{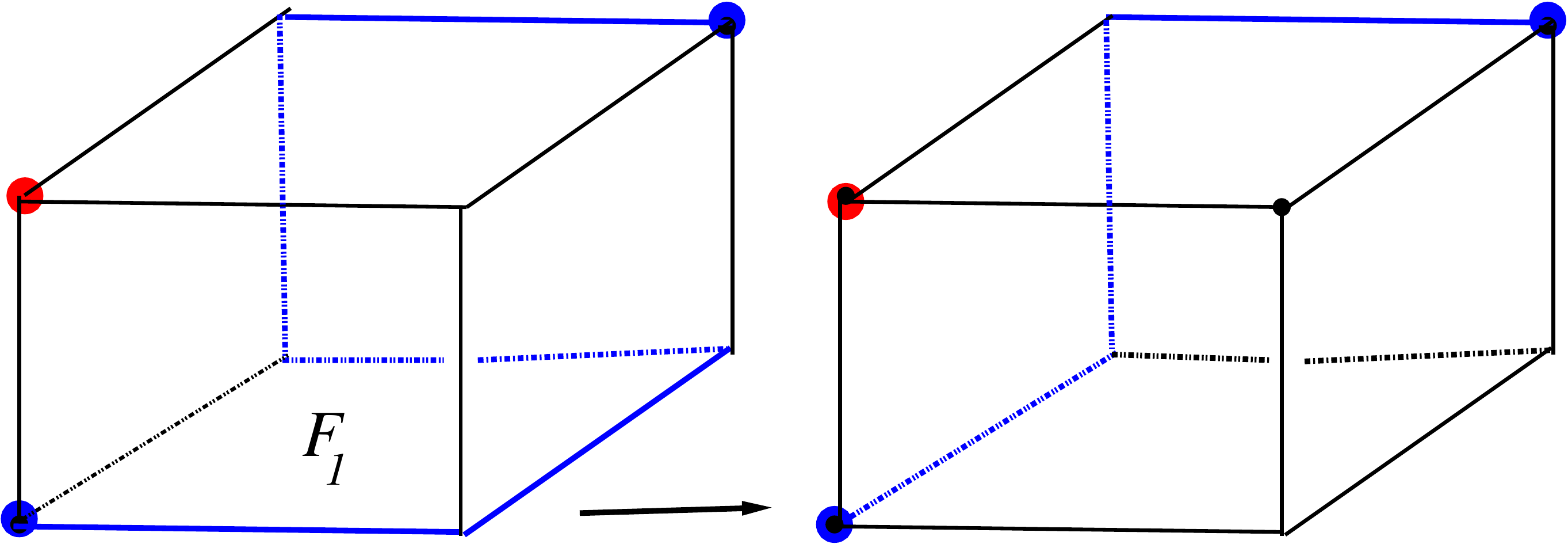}
\end{center}
\caption{\sl $K'\cap Q_j$ consists of a five edges path.} 
\label{I3}
\end{figure} 
\item [$(c)$] $K'\cap Q_j$ consists of a four edges path. If three edges are contained in a face $F_1$ and one edge lies on a neighbor face $F_2$, then
applying an (M2)-move on  $F_1$, we obtain that $K'\cap Q_j$ can be reduce to a two edges path (see Figure \ref{I4}).
\begin{figure}[h] 
 \begin{center}
 \includegraphics[height=3cm]{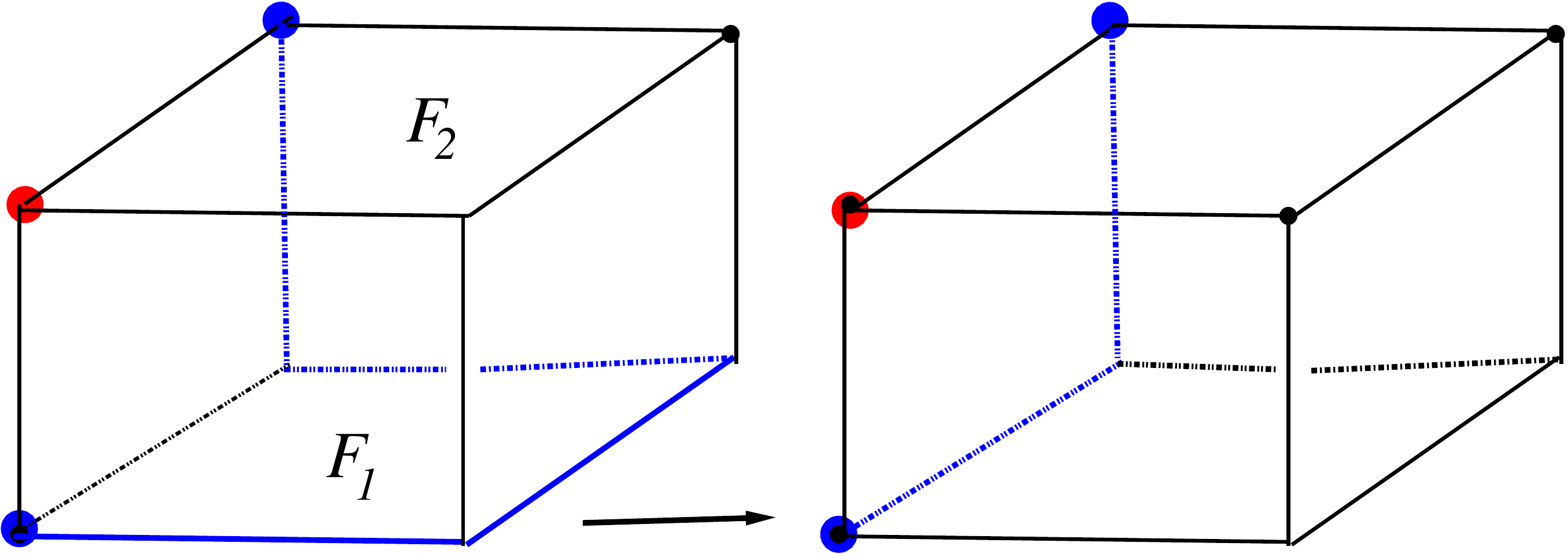}
\end{center}
\caption{\sl $K'\cap Q_j$ consists of a four edges path.} 
\label{I4}
\end{figure} 
Now, if two edges are contained in a face $F_1$ and two edges lie on a neighbor face $F_2$, then
we apply an (M2)-move on  $F_1$ to obtain a new path consisting of four edges. Next we apply again an (M2)-move on  $F_2$ to obtain that $K'\cap Q_j$ 
can be reduced to a two edges path (see Figure \ref{I5}).
\begin{figure}[h] 
 \begin{center}
 \includegraphics[height=3cm]{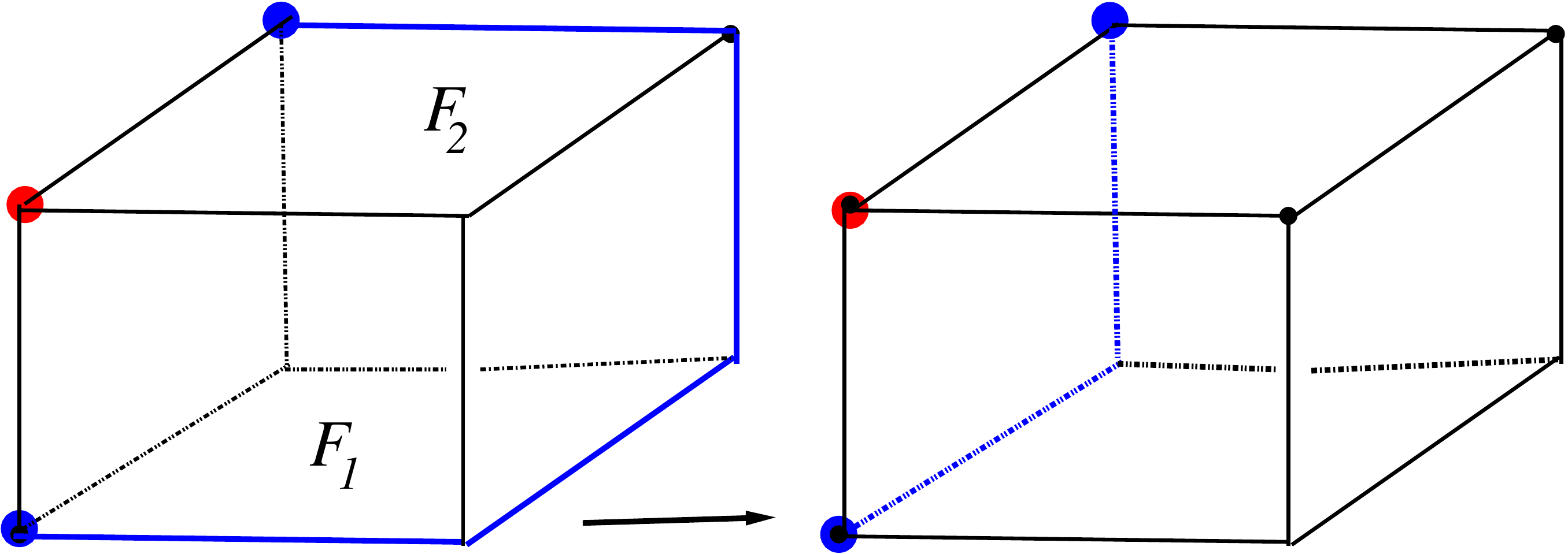}
\end{center}
\caption{\sl $K'\cap Q_j$  consists of a four edges path.} 
\label{I5}
\end{figure} 
\end{description}
\item [Case 2.] Suppose that $K\cap Q_j$ consists of an edge. Then $K'\cap Q_j$ may be contained in either 1 or 2 faces.
\begin{description}
 \item [$(a)$] Three edges of $K'\cap Q_j$ are contained in a face $F_1$ and two edges lie on the neighbor face $F_2$, then
we apply an (M2)-move on  $F_1$ to obtain that $K'\cap Q_j$ can be reduced to a one edge path (see Figure \ref{I6}).
\begin{figure}[h] 
 \begin{center}
 \includegraphics[height=3cm]{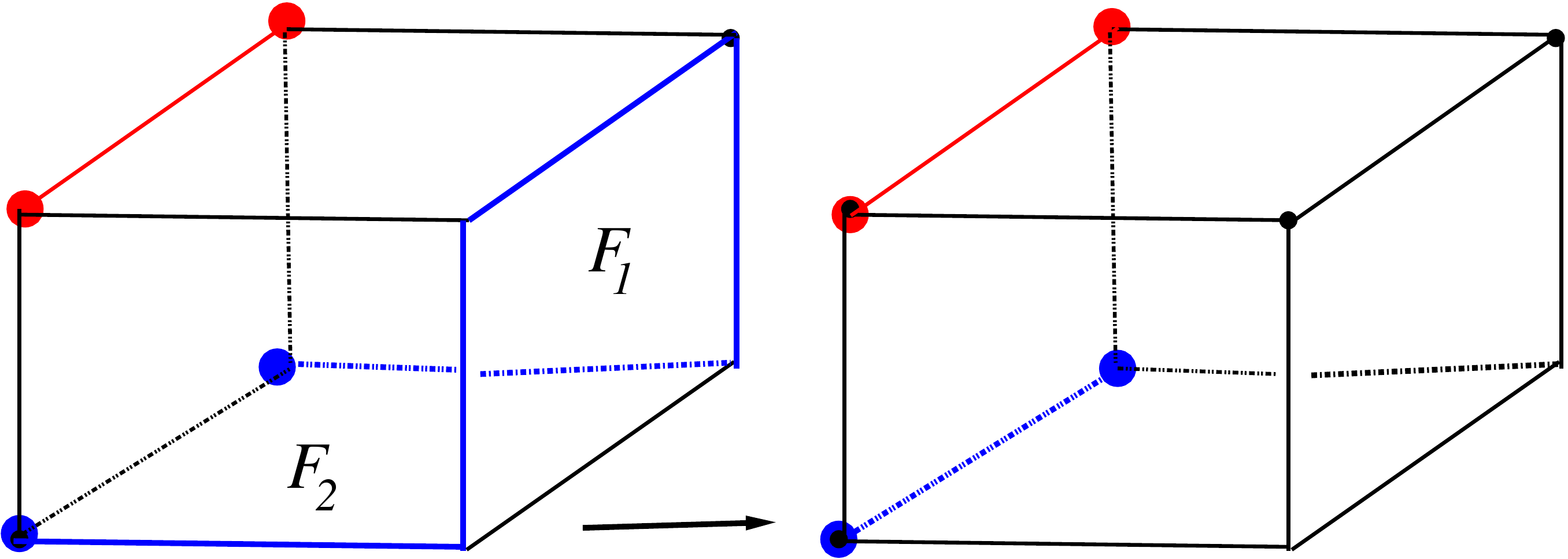}
\end{center}
\caption{\sl $K'\cap Q_j$  consists of a five edges path.} 
\label{I6}
\end{figure} 
\item [$(b)$] Three edges of $K'\cap Q_j$ are contained in a face $F_1$ and one edge lies on the neighbor face $F_2$. Then
we apply an (M2)-move on  $F_1$ to obtain that $K'\cap Q_j$ can be reduced to a two edges path (see Figure \ref{I7}).
\begin{figure}[h] 
 \begin{center}
 \includegraphics[height=3cm]{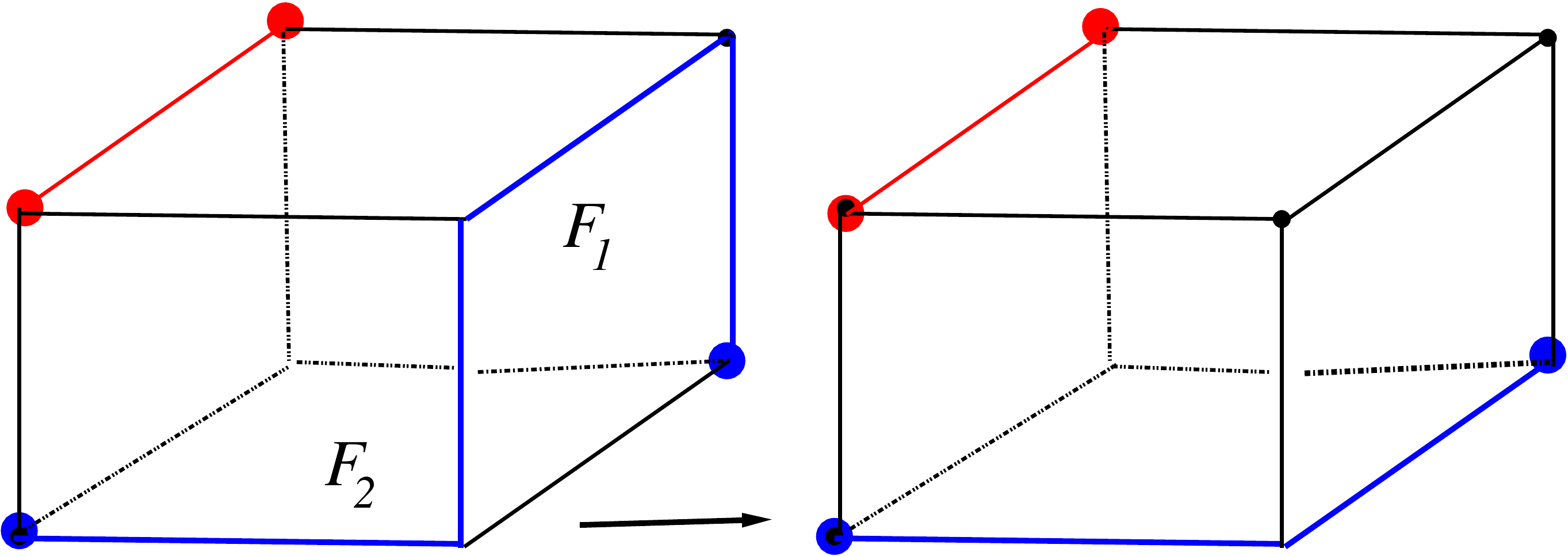}
\end{center}
\caption{\sl $K'\cap Q_j$  consists of a four edges path.} 
\label{I7}
\end{figure} 
\end{description}
\item [Case 3.] Suppose that $K\cap Q_j$ consists of two neighboring edges. In this case, $K'\cap Q_j$ is contained in one face, hence $K'\cap Q_j$ 
is equal to a path consisting of at most three edges. Then, applying an (M2)-move if necessary, we can assume that $K'\cap Q_j$ consists of either one or two edges. 
\end{description}

\noindent This proves claim 2. $\square$\\

By claim 2, we will assume that $K'\cap Q_j$ is equal to a path consisting of at most three edges. Remember that by hypothesis, $K\cap Q_j$
consists of a path of at most two neighboring edges. Next, we will construct $F_j$ considering
 all possible cases of both $K\cap Q_j$ and $K'\cap Q_j$.

\begin{description}
\item  [Case 1.] Suppose that $K\cap Q_j$ consists of a vertex $v$. By claims 1 and 2, we have that $K'\cap Q_j$ consists of a path of either two or three edges.
\begin{description}
\item [$(a)$] $K'\cap Q_j$ is a three edges path which can not be reduced by (M2)-move (see
Figure \ref{f1}). In this case, up to a face boundary move, we may assume that $v$ and two neighboring edges 
of $K'\cap Q_j$ lie on the same face $F$. Then $F_j=F$. Notice that the remaining edge of $K'\cap Q_j$ is contained in a neighbor cube $Q_t$ of $Q_j$ 
which must belong to $B$, hence this edge will be lie on the corresponding $F_t$.
\begin{figure}[h] 
 \begin{center}
 \includegraphics[height=3cm]{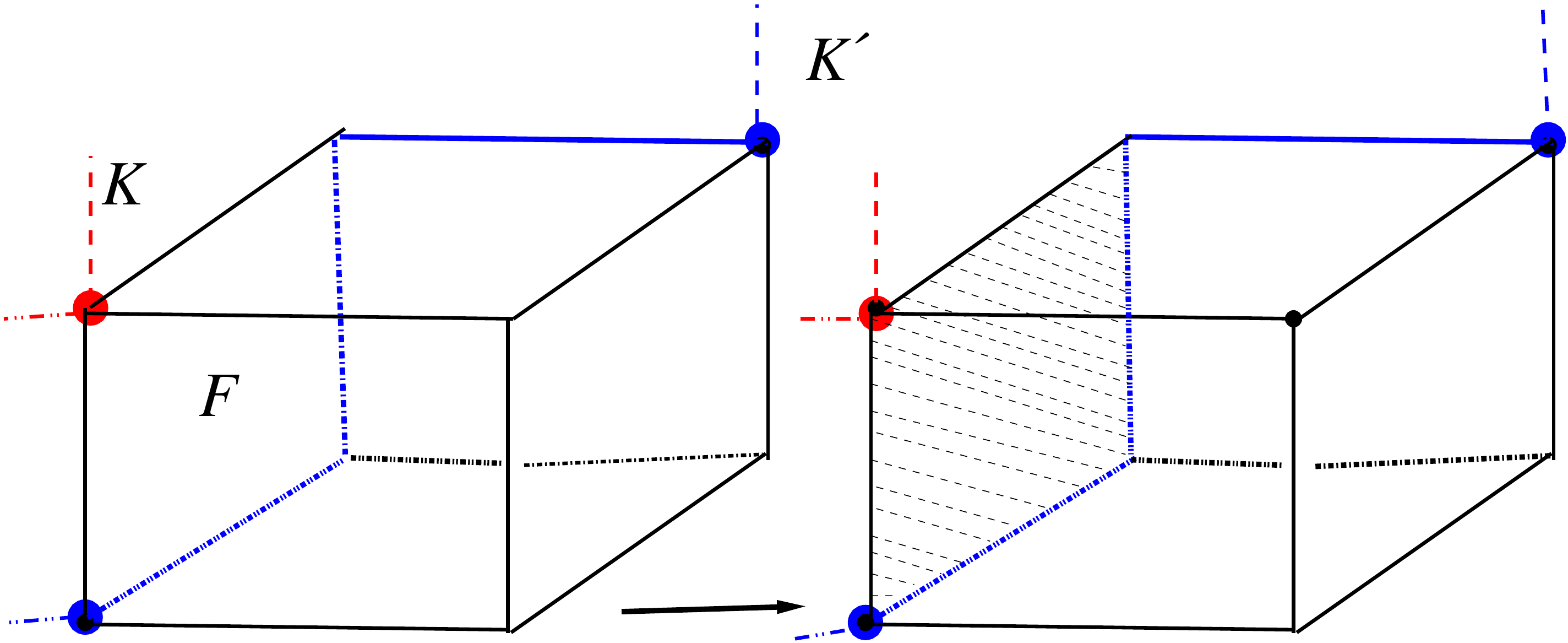}
\end{center}
\caption{\sl $K\cap Q_j$  is a vertex and $K'\cap Q_j$ is a three edges path.} 
\label{f1}
\end{figure} 

\item [$(b)$] $K'\cap Q_j$ is a two edges path. We have the following possible cases. 
\begin{enumerate}
 \item Both $K'\cap Q_j$ and $v$ lie on the same face $F\subset Q_j$. Then $F_j=F$. See Figure \ref{f2}.
\begin{figure}[h] 
 \begin{center}
 \includegraphics[height=3cm]{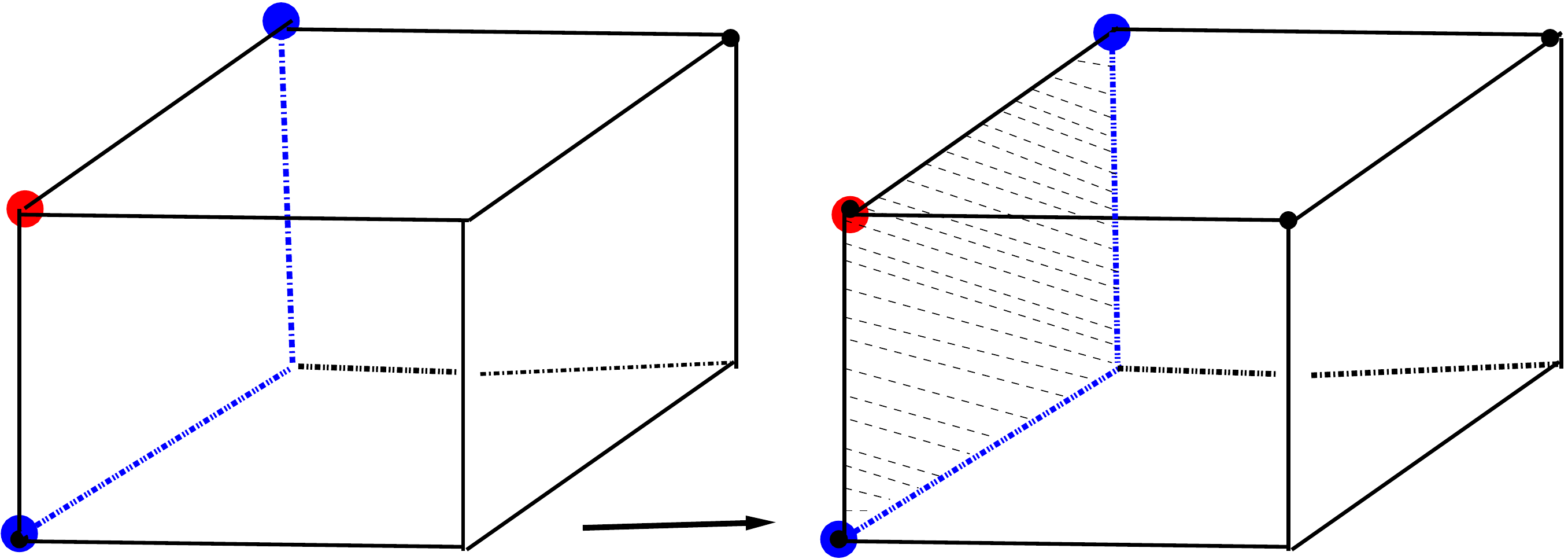}
\end{center}
\caption{\sl Both $K'\cap Q_j$ and $v$ lie on the same face.} 
\label{f2}
\end{figure}
\item $K'\cap Q_j$ and $v$ lie on opposite faces. We have two possibilities which are equivalent via an (M2)-move (see Figure \ref{f3}).
\begin{figure}[h] 
 \begin{center}
 \includegraphics[height=3cm]{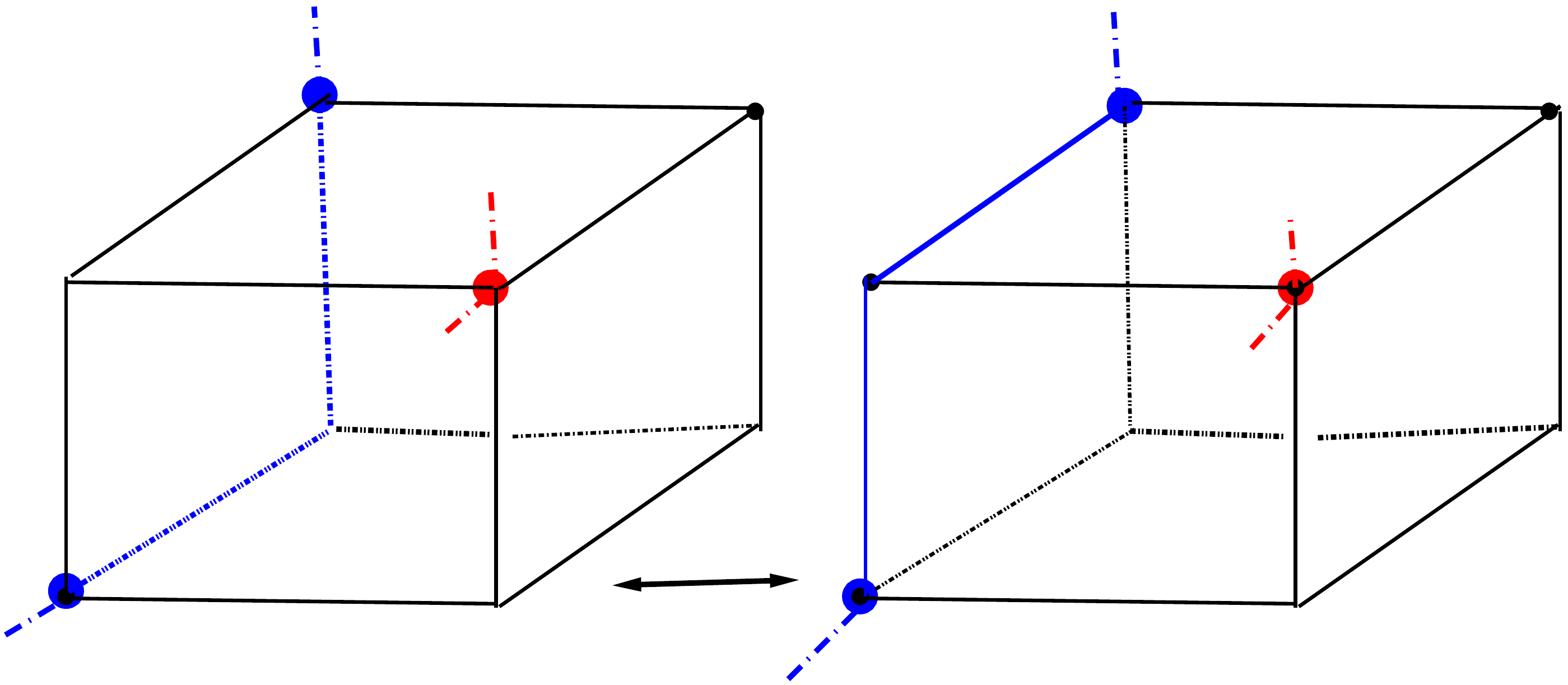}
\end{center}
\caption{\sl $K\cap Q_j$ and $K'\cap Q_j$ lie on opposite faces.} 
\label{f3}
\end{figure}
Now, we have that there exist neighbor cubes of $Q_j$ such that they intersect both $K$ and $K'$. Notice that this can only happens in the
configuration described on Figure \ref{f4}. Then applying  (M2)-moves on these neighbor cubes, we have that $Q_j$ is superfluous, hence $F_j=\emptyset$.
\begin{figure}[h] 
 \begin{center}
 \includegraphics[height=3cm]{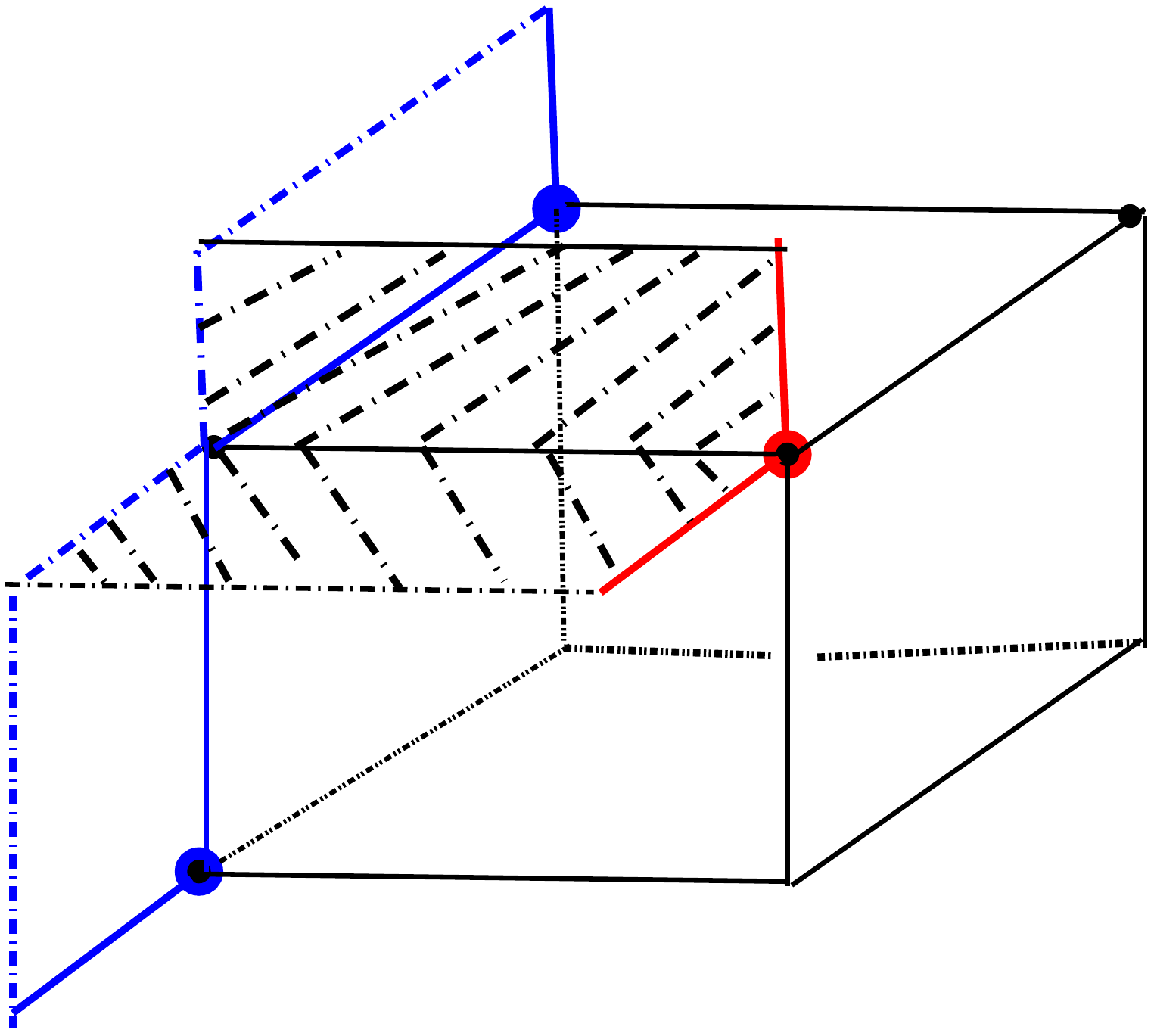}
\end{center}
\caption{\sl $F_j$ is empty.} 
\label{f4}
\end{figure}
\item The last case appears in Figure \ref{f5}. By the same argument of claim 1 (b), this configuration is not possible.
\begin{figure}[h] 
 \begin{center}
 \includegraphics[height=3cm]{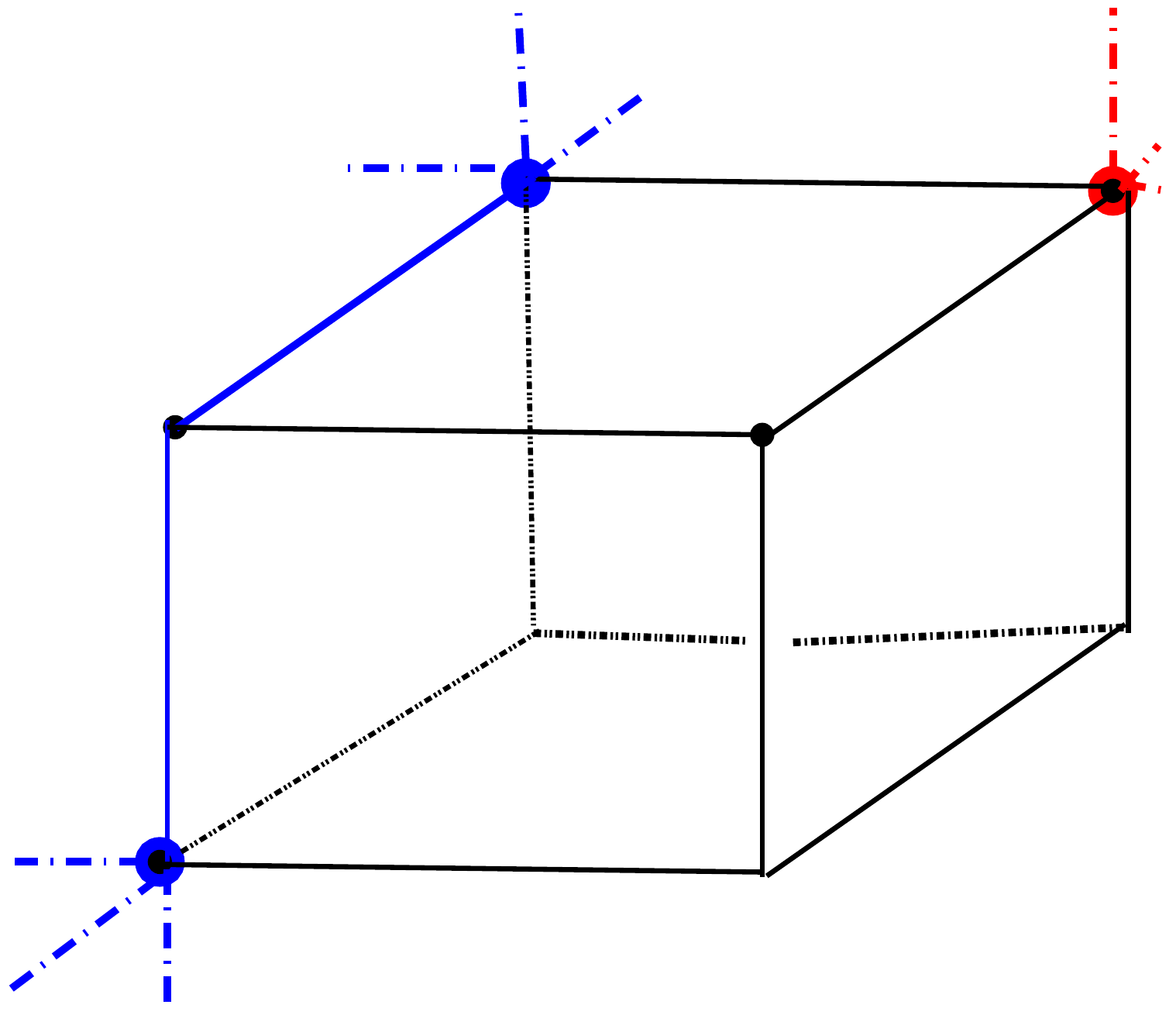}
\end{center}
\caption{\sl This configuration is not possible.} 
\label{f5}
\end{figure}
\end{enumerate}

\end{description}

\item  [Case 2.] Suppose that $K\cap Q_j$ consists of one edge. By claim 2, we can assume that $K'\cap Q_j$ consists of a path of either one, two or three edges.
\begin{description}
\item [$(a)$] $K'\cap Q_j$ is a three edges path which can not be reduced by an (M2)-move to a path consisting of either one or two edges. 
Since $K\cap K'=\emptyset$, then  $K'$ must be contained in two neighboring faces 
of $Q_j$ which intersect $K$ (see Figure \ref{F}). So $F_j$ is the union of these two neighboring faces. 
\begin{figure}[h] 
 \begin{center}
 \includegraphics[height=3cm]{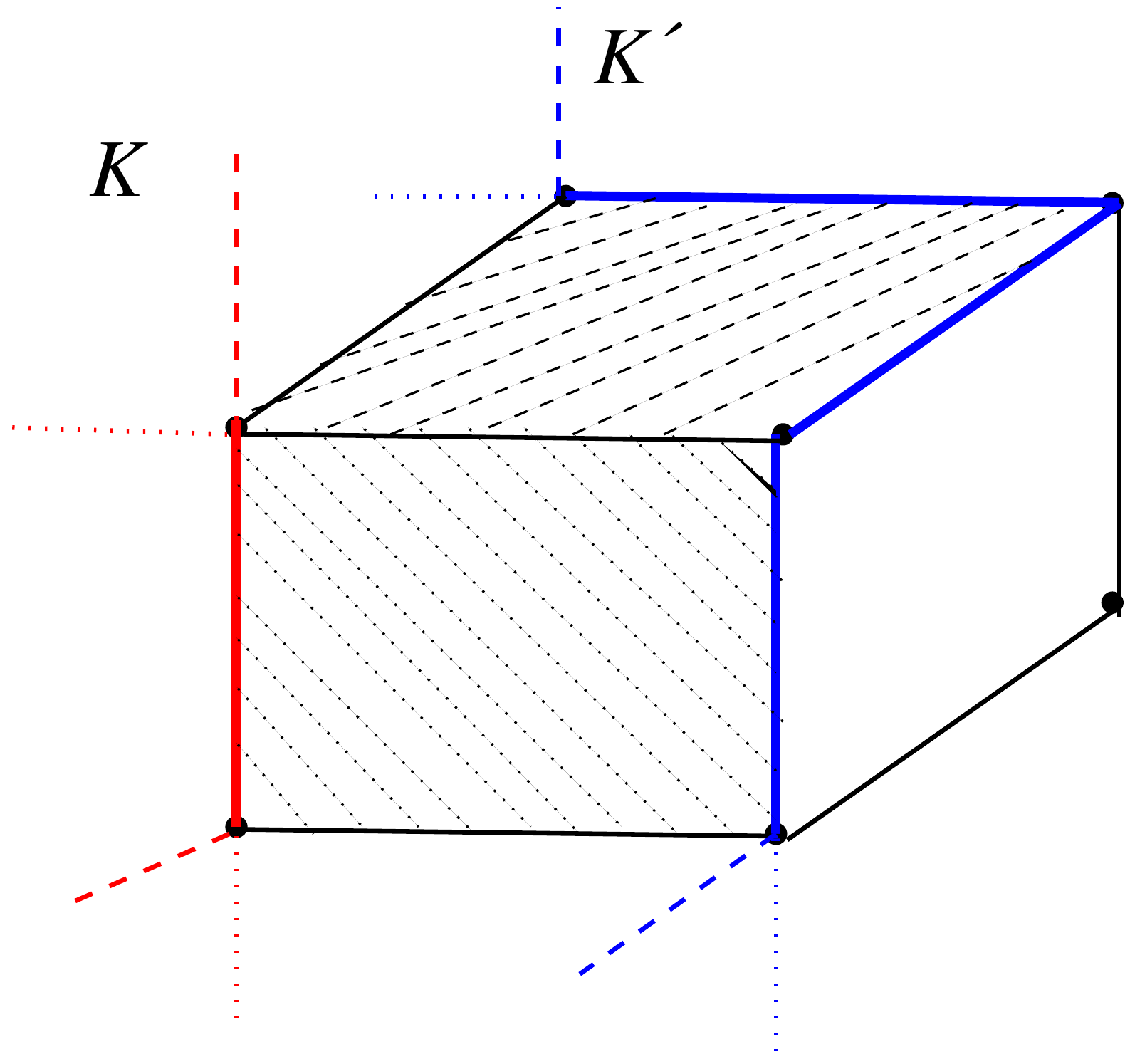}
\end{center}
\caption{\sl $K'\cap Q_j$ consists of a three edges path.} 
\label{F}
\end{figure} 
 \item [$(b)$] $K'\cap Q_j$ is a two edges path.
\begin{enumerate}
\item Two edges of $K\cup K'$ lie on a face $F\subset Q_j$. Then $F_j=F$ (see Figure \ref{F2}). Notice that the remaining edge
 is contained in a neighbor cube $Q_t$ of $Q_j$ 
which must belong to $B$, hence this edge will be lie on the corresponding $F_t$.
\begin{figure}[h] 
 \begin{center}
 \includegraphics[height=3cm]{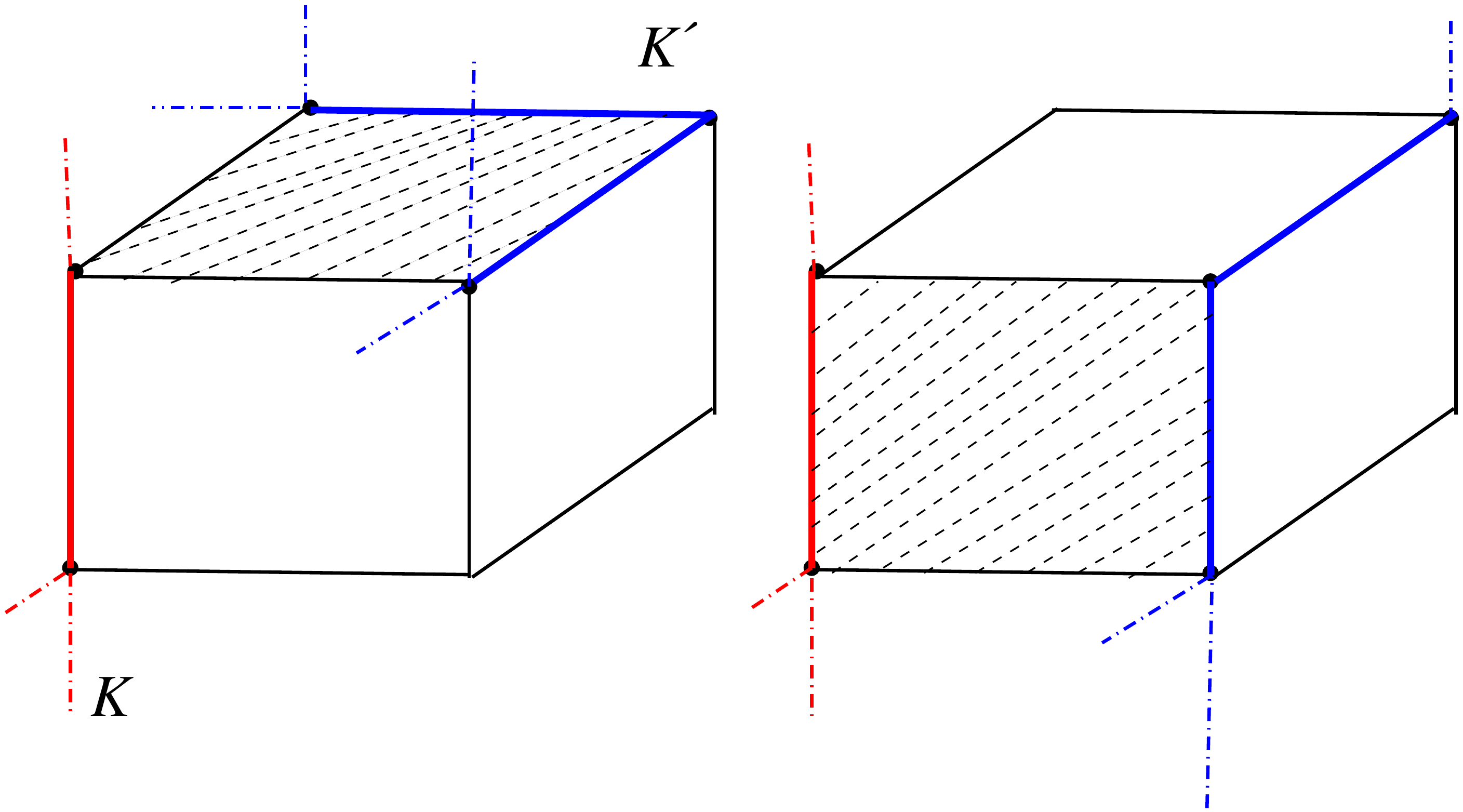}
\end{center}
\caption{\sl $K\cap Q_j$ consists of an edge.} 
\label{F2}
\end{figure} 
\item The last possible configuration appears in Figure \ref{f6}. Then applying an (M2)-move if necessary, 
it is equivalent to the configuration described above (see Figure \ref{F2}).
\begin{figure}[h] 
 \begin{center}
 \includegraphics[height=3cm]{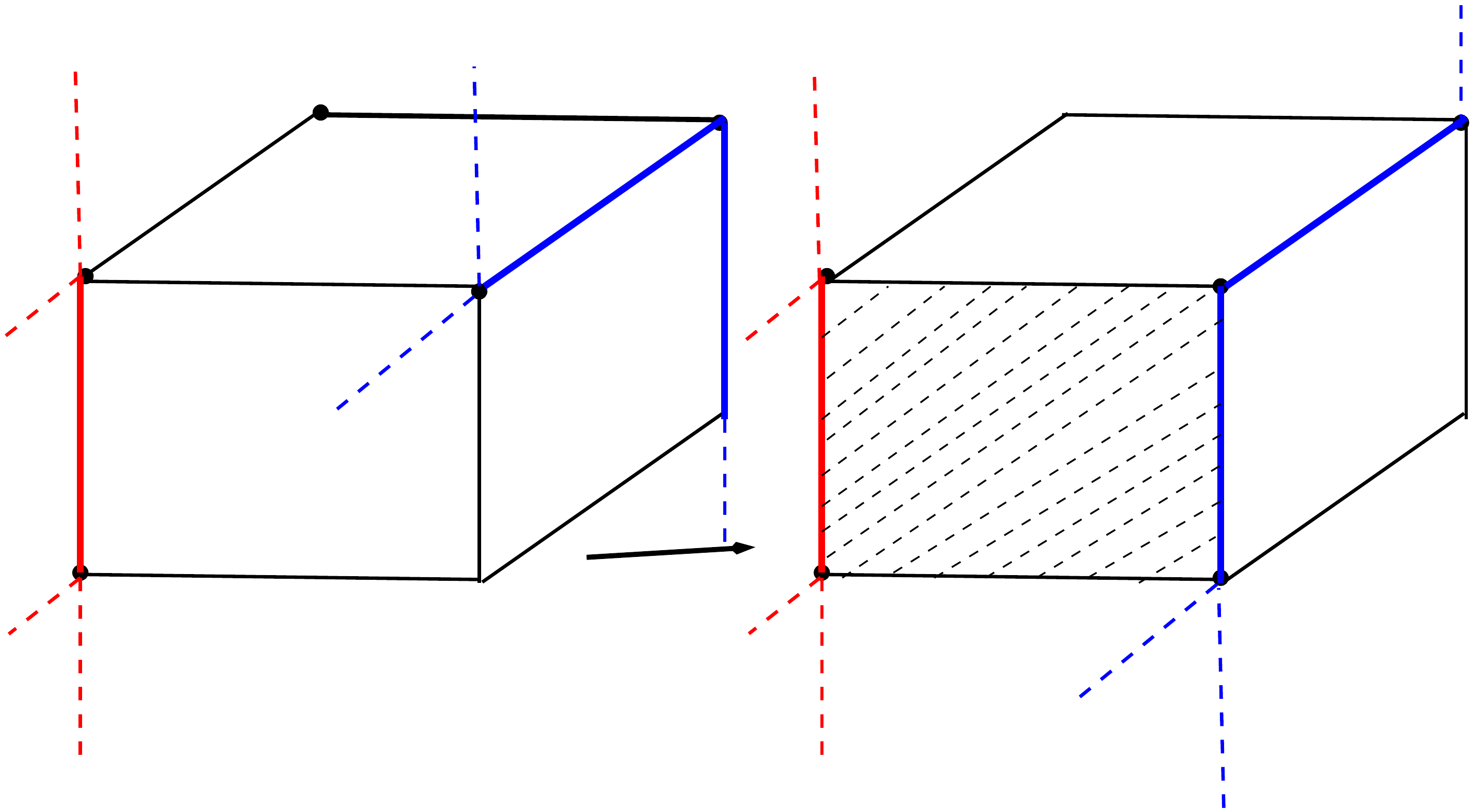}
\end{center}
\caption{\sl Equivalent configuration.} 
\label{f6}
\end{figure} 
\end{enumerate}

\item [$(c)$] $K'\cap Q_j$ is an one edge path.
\begin{enumerate}
 \item Both $K'\cap Q_j$ and $K\cap Q_j$ lie on the same face $F\subset Q_j$. Then $F_j=F$ (see Figure \ref{f7}).
\begin{figure}[h] 
 \begin{center}
 \includegraphics[height=3cm]{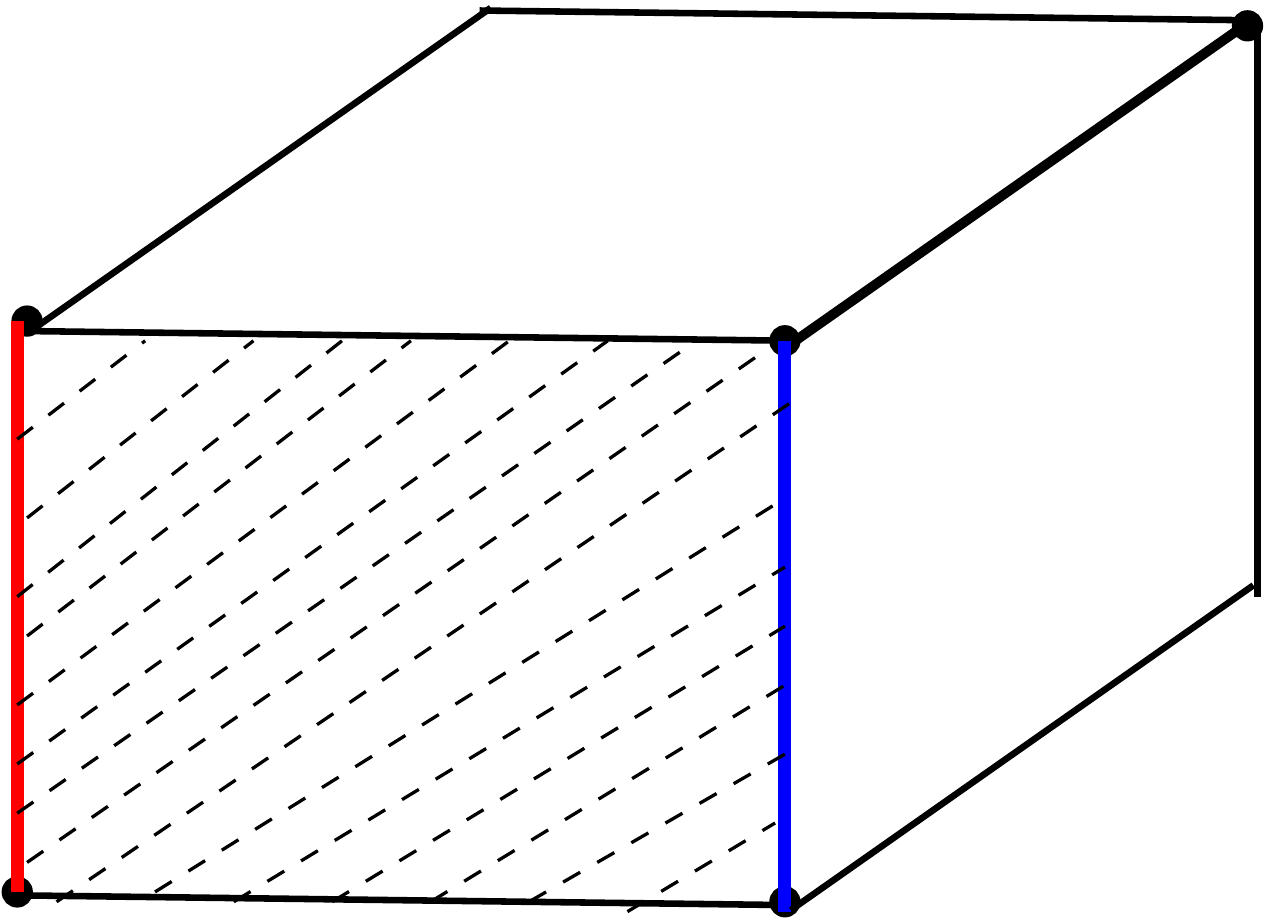}
\end{center}
\caption{\sl Both $K'\cap Q_j$ and $K\cap Q_j$ lie on the same face.} 
\label{f7}
\end{figure} 

\item $K'\cap Q_j$ and $K\cap Q_j$ are two opposite edges of $Q_j$, \emph{i.e.}, consider the vertices $w_1$, $w_2\in K\cap Q_j$ and the vertices 
$v_1$, $v_2\in K'\cap Q_j$, then
$v_1$, $w_1$ lie on the face $F_1\subset Q_j$ and
$v_2$, $w_2$ lie on the opposite face $F_2\subset Q_j$; such that $v_1$, $w_1$ are opposite vertices on $F_1$, and $v_2$, $w_2$ are opposite vertices on $F_2$. 
See Figure \ref{F5}. In other words, the edge $\overline{v_1v_2}$ is opposite to the edge $\overline{w_1w_2}$.
\begin{figure}[h] 
 \begin{center}
 \includegraphics[height=3.5cm]{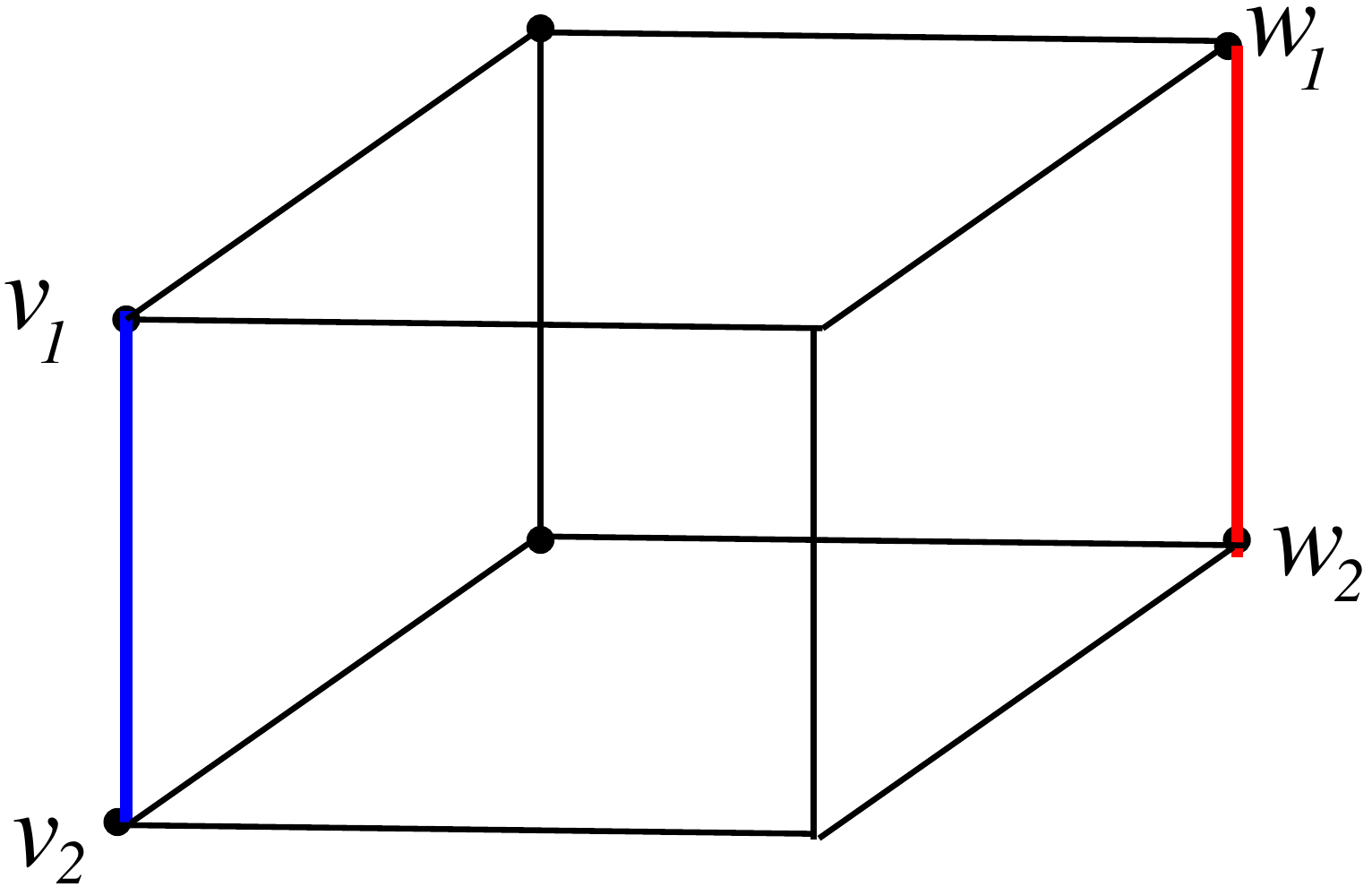}
\end{center}
\caption{\sl $K'\cap Q_j$ and $K\cap Q_j$ are opposite edges.} 
\label{F5}
\end{figure}

Observe that there exists a finite sequence of neighbor cubes, $Q_{i+1},\ldots,Q_j,\ldots,Q_{r-1}$, keeping this configuration. However, at the cubes $Q_i$ and $Q_r$
this configuration changes. This implies that each intersection $K'\cap Q_i$, $K'\cap Q_r$ 
consists of two neighboring edges.
We can assume, applying induction if necessary, that this happens in the neighbor cubes of $Q_j$; {\it i.e.} $i=j-1$ and $r=j+1$. Let
$F_1=Q_{j-1}\cap Q_j$ and $F_2=Q_{j}\cap Q_{j+1}$. 

Consider the cube $Q_{j-1}$. Then we have two possible configurations. 
\begin{enumerate}
\item Let $w_1$, $w_3$, $v_1$, $v_3$, $v\in Q_{j-1}$ be vertices  satisfying the following.
$w_1$, $w_3\in K$, $v_1$, $v_3$, $v\in K'$, the edge $\overline{v_1v_3}$ is opposite to the edge $\overline{w_1w_3}$
and the vertices $v_3$, $v$, $w_3$ lie on the face $F'_1$ opposite to  $F_1$ (see Figure \ref{F6}).\\

\begin{figure}[h]  
 \begin{center}
 \includegraphics[height=4.5cm]{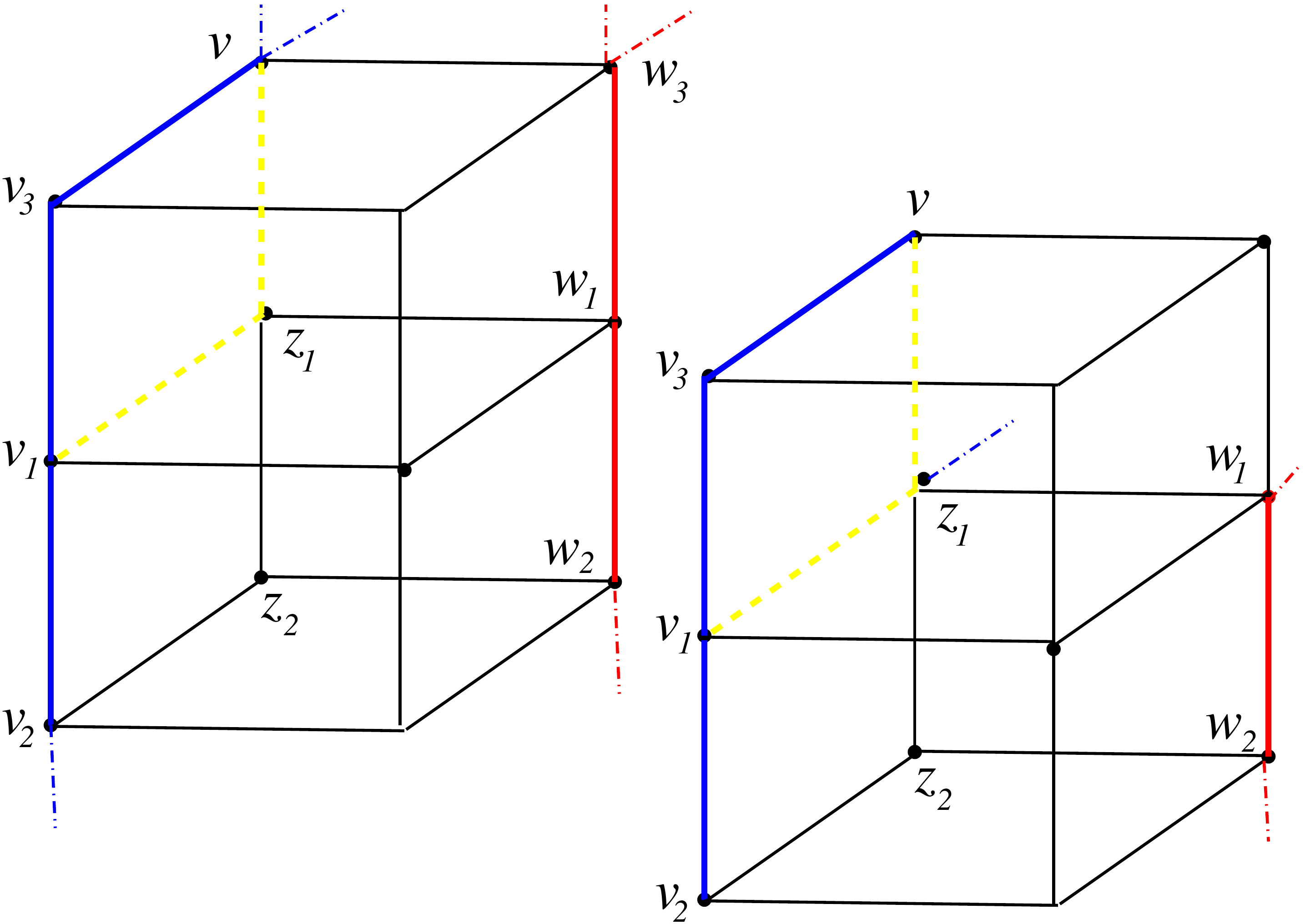}
\end{center}
\caption{\sl Possible arrangements.} 
\label{F6}
\end{figure}
Notice that we have the same configuration described on case $2(b)1$. Applying an (M2)-move on $Q_{j-1}$, we get that the edges 
$\overline{vv_3}\cup\overline{v_1v_3}$ are replaced by the edges $\overline{vz_1}\cup\overline{v_1z_1}$, hence  $F_{j-1}$ is exactly the face
of $Q_{j-1}$ containing $\overline{vz_1}\cup\overline{w_3w_1}$. Now, consider the cube $Q_j$, so it has the same configuration analyzed on
 $2(b)1$. Then, we apply an (M2)-move to obtain $F_j$. 

\item Suppose that $K\cap Q_{j-1}=w_1$. Then we apply an (M2)-move on $Q_{j-1}$ (see case$2(b)1$), so the edges 
$\overline{vv_3}\cup\overline{v_1v_3}$ are replaced by the edges $\overline{vz_1}\cup\overline{v_1z_1}$. Now, there exists another cube $Q\in B$ such that
$v$, $z_1$, $w_1\in Q\cap Q_{j-1}$; hence $Q_{j-1}$ is superfluous (see Claim 1). Observe that $Q_j$ has the same configuration described on $2(b)1$, then
applying an (M2)-move, we obtain $F_j$. Notice that $F_{j-1}=\emptyset$.
\end{enumerate}
\item The last possible configuration appears in Figure \ref{f8}. As we can see, it is very similar to the case 1$(b)2$. Then, by the same 
argument $F_j=\emptyset$.
\begin{figure}[h] 
 \begin{center}
 \includegraphics[height=3cm]{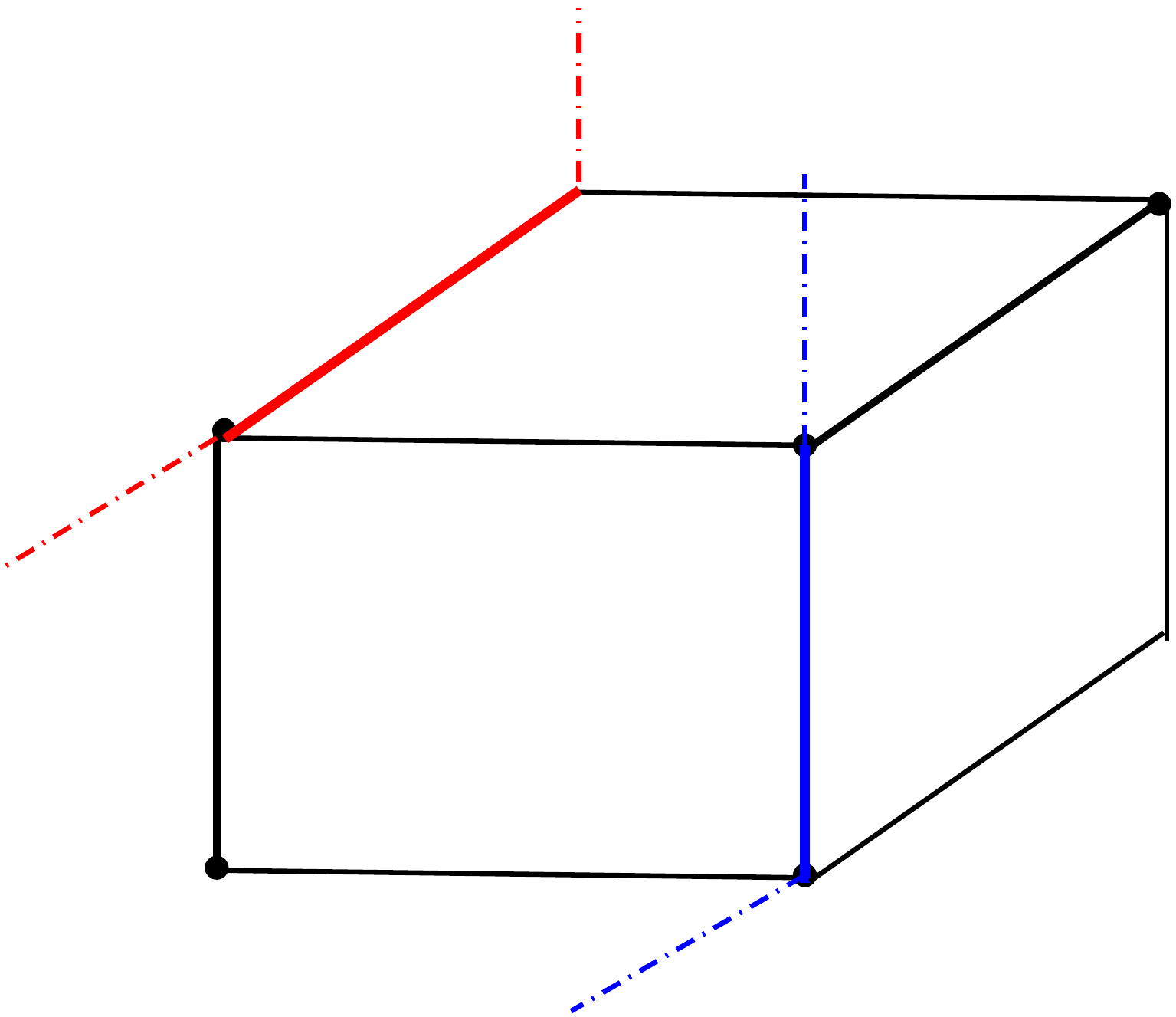}
\end{center}
\caption{\sl $F_j=\emptyset$.} 
\label{f8}
\end{figure} 
\end{enumerate}
\end{description}
\item  [Case 3.] Suppose that $K\cap Q_j$ consists of two neighboring edges. Then $K\cap Q_j$ lies on a face $F\subset Q_j$.
\begin{description}
\item [$(a)$] $K'\cap Q_j$ is a three edges path. Since $K\cap K'=\emptyset$, then $K'$ lies on five possible edges, four of them
belong to the opposite face $F'$ of $F$,  and the remaining edge $e$. Hence $K'$ consists of $e$ and two edges $l_1,\,l_2\subset F'$ (see Figure \ref{f9}). 
Applying an (M2)-move if necessary, we can assume that $e$ and $l_1$ belong to a face $F_1$ and $l_2$ and an edge of $K\cap Q_j$ are contained in a 
neighbor face $F_2$. Notice that the remaining edge of $K\cap Q_j$ is contained in a neighbor cube of $Q_j$ which belongs to $B$, where it is considered.
Then $F_j=F_1\cup F_2$.  
\begin{figure}[h] 
 \begin{center}
 \includegraphics[height=3cm]{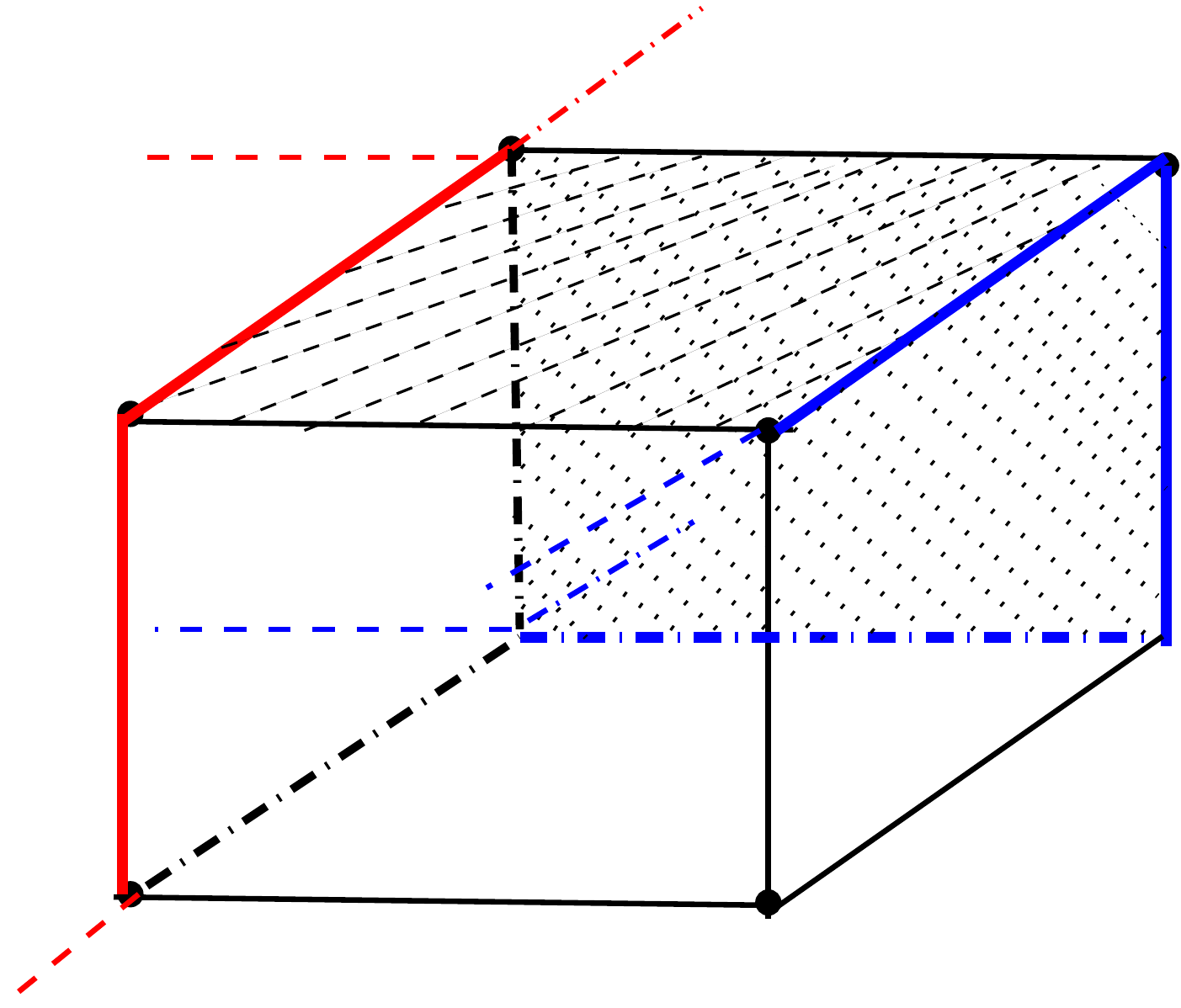}
\end{center}
\caption{\sl $K\cap Q_j$ and $K'\cap Q_j$ lie on two neighboring faces.} 
\label{f9}
\end{figure} 

\item [$(b)$] $K'\cap Q_j$ is a two edges path.
\begin{enumerate}
 \item $K\cap Q_j$ and $K'\cap Q_j$ lie on two neighboring faces (faces sharing a common edge) of $Q_j$. 
Let $F_j$ be the union of these two faces. See Figures \ref{F1} and \ref{f10}.
\begin{figure}[h] 
 \begin{center}
 \includegraphics[height=3cm]{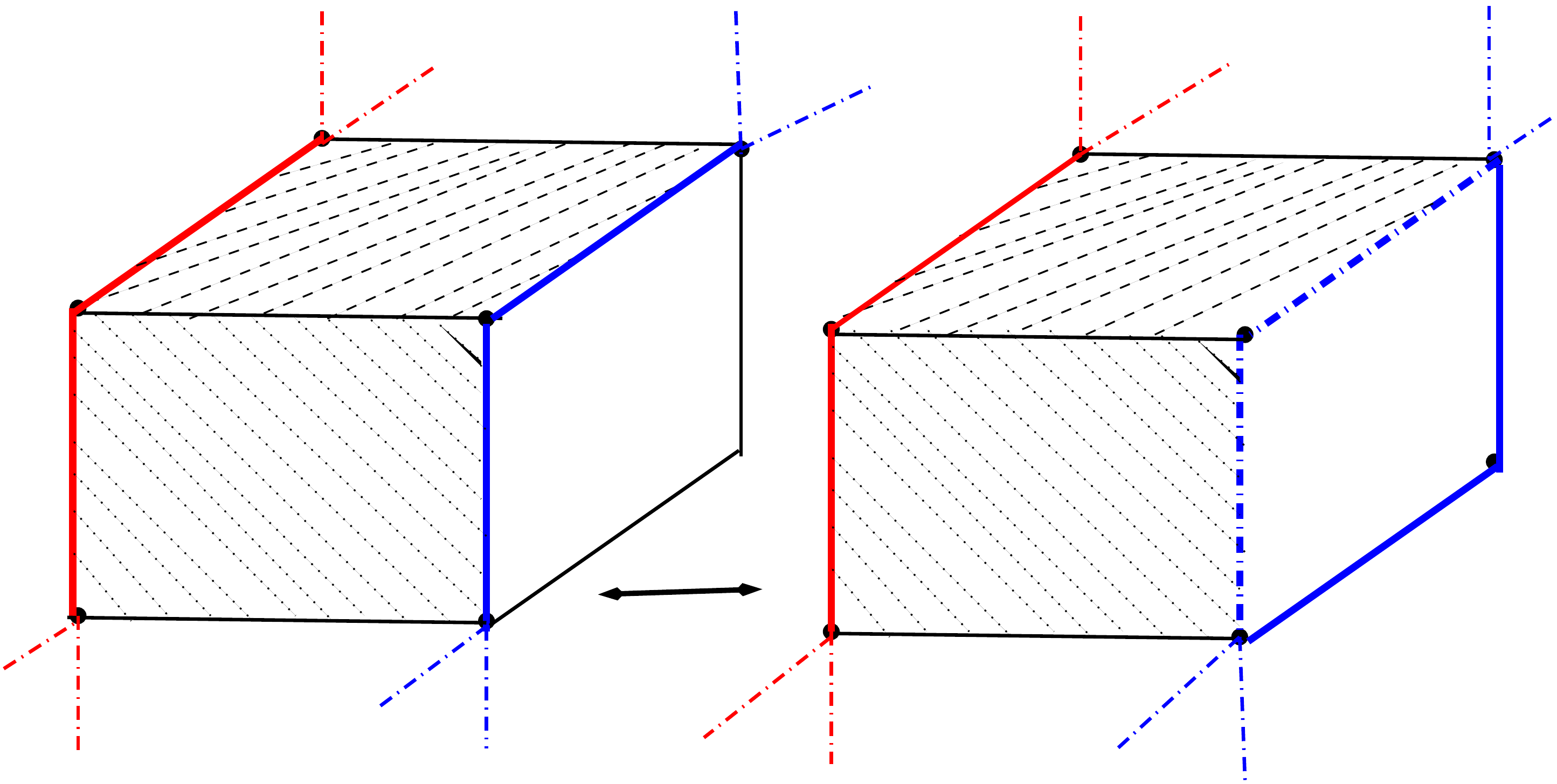}
\end{center}
\caption{\sl $K\cap Q_j$ and $K'\cap Q_j$ lie on two neighboring faces.} 
\label{F1}
\end{figure} 

\begin{figure}[h] 
 \begin{center}
 \includegraphics[height=3cm]{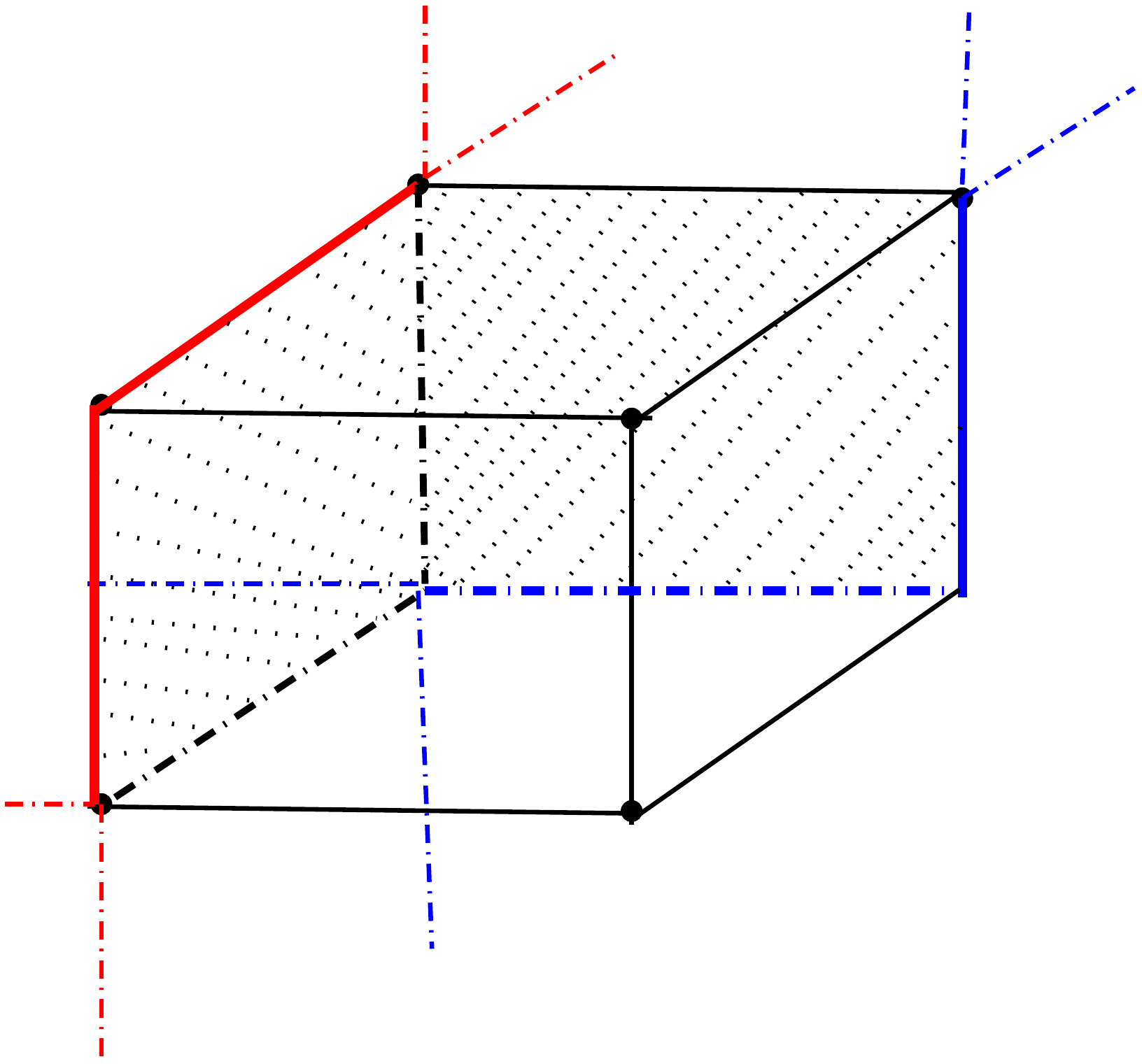}
\end{center}
\caption{\sl $K\cap Q_j$ and $K'\cap Q_j$ lie on two neighboring faces.} 
\label{f10}
\end{figure} 

\item $K'\cap Q_j$ lies on the opposite face $F'$ of $F$, in such a way that there is only one face $F_1\subset Q_j$ such that
both $F_1\cap K$ and $F_1\cap K'$ consists on an edge, respectively. Notice that the remaining edges of $K$ and $K'$ belong to
neighboring cubes of $Q_j$, where they are considered. Then $F_j=F_1$. See Figure \ref{f11}
\begin{figure}[h] 
 \begin{center}
 \includegraphics[height=3cm]{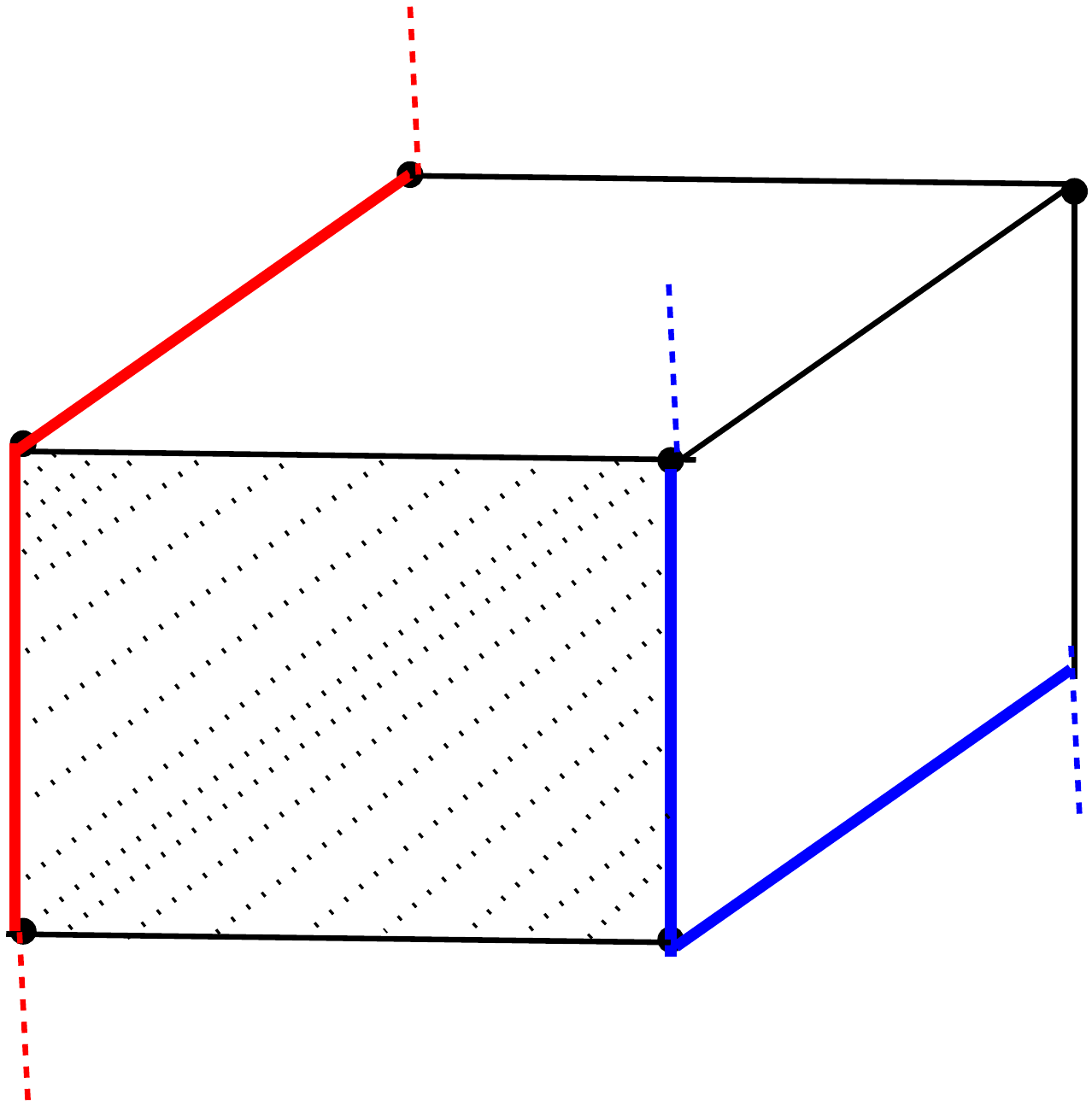}
\end{center}
\caption{\sl $K\cap Q_j$ and $K'\cap Q_j$ lie on opposite faces.} 
\label{f11}
\end{figure} 
\end{enumerate}

\end{description}

\end{description}

Let  ${\cal {F}}=\cup{F_j}$. By the construction, $F_j$ consists of either one face or two neighboring faces which are contained in the  2-skeleton of ${\cal{N}}_k$,
hence ${\cal{F}}$ is contained in the  2-skeleton of ${\cal{N}}_k$.\\

As a consequence of the previous analysis, we know that if we apply an (M2)-move on a cube $Q_j$, then this movement does not alter the configuration on its neighbor cubes; in other words,
the choice of $F_j$ depends only of the configuration given on $Q_j$. Now the boundary of each $F_j$ is composed by  edges belonging to  $K$ and $K'$  and 
two disjoint edges $e_j$ and $e'_j$ that do not belong to
neither $K$ nor $K'$. Thus, if we take the cubes $Q_j$, $Q_{j+1}$ and we construct $F_j$ and $F_{j+1}$, then 
the intersection $F_j\cap F_{j+1}$ is one of these two edges $e_j$, $e'_j$; moreover, since $K$ and $K'$ are closed, connected 1-manifolds, then
each of these two edges belong exactly to two cubes. 
This implies that ${\cal{F}}$ is a surface whose boundary
consists of two connected components, namely $K$ and $K'$.\\

Now, we will carry the knot $K'$ onto the knot $K$ via a finite number of cubulated moves. Notice that ${\cal{F}}$ is the union 
of $m$ squares $F_1,\,\ldots F_m$ which are enumerated according to the numbering provided by the set $B$. In particular, the square $F_{n+1}$ is a neighbor square of $F_n$ for all $n$, 
\emph{i.e.}, $F_{n+1}\cap F_n$ consists of an edge $l_n$.\\

We will use induction on $m$. Consider $F_1$. We apply an (M2)-move on $F_1$ in such a way that the edges belonging to $K'$ replace the others
(by construction $F_1\cap K\neq\emptyset$ and $F_1\cap K'\neq \emptyset$).
Next, we consider $F_2$. By the previous step, $F_1$ and $F_2$ share an edge belonging to $K'$. Then we apply again an (M2)-move, so
the edges belonging to $K'$ replace the others. We continue this finite process in thi way. Notice that if $l\subset  K$  
then $l\subset \partial {\cal{F}}$, thus $l$ is not a common edge of two squares $F_i$ and $F_j$ belonging to ${\cal{F}}$, hence if 
$l$ is replaced by an edge belonging to $K'$, then this replacement will be kept in the following steps. Therefore, the result follows. $\square$

\subsubsection{Smooth case}

Suppose that  $\widetilde{K}_{1}$ and $\widetilde{K}_{2}$ are isotopic  smooth 1-knots in  $\mathbb{R}^{3}$. 
Let $J$ be the isotopy cylinder of $\widetilde{K}_{1}$ and $\widetilde{K}_{2}$.
Then $J$ is a smooth submanifold of codimension two in $\mathbb{R}^{4}$. 
By Theorem \ref{trace}, there exists an isotopic copy of $J$, say $J'$, contained in the $2$-skeleton of the canonical 
cubulation $\cal{C}$ of $\mathbb{R}^{4}$, such that $J'$ is sliced by connected level sets of $p$. Moreover there exist integer numbers 
$m_1$ and $m_2$ and cubic knots $K'_1$ and $K'_2$ isotopic to $\widetilde{K}_{1}$ and $\widetilde{K}_{2}$, respectively;
such that $p^{-1}(t)\cap J' = K'_1$ for all $t\leq m_1$ and
$p^{-1}(t)\cap J' = K'_2$ for all $t\geq m_2$.\\

Our goal it to prove that  $K'_1$ is equivalent to $K'_2$ by cubulated moves,\\
  
Let $p:\mathbb{R}^{4}\hookrightarrow\mathbb{R}$ be the projection onto the last coordinate. 
Thus $p^{-1}(t)=\mathbb{R}^3_t$ is an affine hyperplane
parallel to the space $\mathbb{R}^3\times\{0\}$.

\begin{lem}
 Each hyperplane $\mathbb{R}^3_t$ has a canonical cubulation. 
\end{lem}

\noindent{\it Proof.} Since  $\mathbb{R}^3_t$ is an affine hyperplane parallel to the
hyperplane ${\cal{P}}=\mathbb{R}^3\times\{0\}\subset\mathbb{R}^4$ and  ${\cal{P}}$ has a canonical cubulation ${\cal{C}}_{\cal{P}}$ given by the restriction of
the cubulation $\cal C$ of $\mathbb{R}^{4}$ to it; {\it i.e.}, ${\cal{C}}_{\cal{P}}$ is the decomposition into cubes which are the images of the unit cube
$I^{3}=\{(x_{1},x_2,x_{3},0)\,|\,0\leq x_{i}\leq 1\}$ by translations by vectors with integer coefficients whose last coordinate is zero. Then
$\mathbb{R}^3_t$ has a canonical cubulation coming from the translation of the cubulation ${\cal{C}}_{\cal{P}}$ by the vector $(0,0,0,t)$. $\square$

\begin{definition}
 A $2$-cell ($2$-face) $F$ of the canonical cubulation $\cal{C}$ of $\mathbb{R}^4$  is called horizontal if $p(F)$ is a fixed integer number. Otherwise
$F$ is called vertical.
\end{definition}

\begin{definition}
 Let $\Sigma$ be a cubulated surface (a surface contained in the $2$-skeleton of $\cal{C}$) and  $P$ be a hyperplane in $\mathbb{R}^4$. We say that 
$P$  and $\Sigma$ intersect transversally, denoted by $P\pitchfork\Sigma$, if $P\cap\Sigma$ is a polygonal curve.
\end{definition}

\begin{lem}
 Let $\mathbb{R}^3_t$, $t\not\in\mathbb{Z}$ be an affine hyperplane. Then  $\mathbb{R}^3_t$ intersects $J'$ transversally.
\end{lem}

\noindent{\it Proof.} Let $\mathbb{P}= \mathbb{R}^3_t$, $t\not\in\mathbb{Z}$. By Theorem \ref{trace},  
$\mathbb{P}\cap J'$ is connected. Let $x\in \mathbb{P}\cap J'$. Then
$x\in F_i$, where $F_i$ is a 2-face of the cubulation $\cal{C}$. Notice that $F_i$ is a vertical $2$-face, since $t\not\in\mathbb{Z}$. 
So, we have two possibilities:
either $x\in\mbox{Int}(F_i)$ or $x$ belongs to an edge of $F_i$.
\begin{enumerate}
\item If $x\in\mbox{Int}(F_i)$ then $\mathbb{P}\cap F_i$ is a linear segment parallel to an edge contained in $\mathbb{R}^3\times\{0\}$.
\item If $x$ belongs to an edge of $F_i$, then there exists another vertical $2$-face $F_j$ such that $x\in F_i\cap F_j$. Thus 
$\mathbb{P}\cap F_i$ is a linear segment $l_i$ parallel to an edge,  and by the same argument  $\mathbb{P}\cap F_j$ is also a linear segment $l_j$,
and  $x\in l_i\cup l_j$. 
\end{enumerate}

\noindent Therefore, the result follows. $\square$\\

\begin{coro}
Let $x$ be a real number such that $x\not\in\mathbb{Z}$. Then the set $p^{-1}(x)\cap J'$ is a knot.
\end{coro}

\noindent{\it Proof.} By the above, $p^{-1}(x)\cap J'$ is a polygonal connected curve. $\square$\\

Now, for each $n\in\mathbb{N}$ we define 
$$
K_{-\frac{1}{2}}(n)= p^{-1}(n-\frac{1}{2})\cap J'
$$ 
and 
$$
K_{\frac{1}{2}}(n)= p^{-1}(n+\frac{1}{2})\cap J'.
$$
Observe that $K_{-\frac{1}{2}}(n)$ and $K_{\frac{1}{2}}(n)$ are cubic knots.\\

Let $\cal{Q}(C)$ be the set of squares ($2$-cells) belonging to ${\cal C}$. Consider the spaces
$$
B_{-\frac{1}{2}}(n)=\cup\{F\in{\cal{Q}(C)}\,|\,F \cap K_{-\frac{1}{2}}(n)\neq\emptyset\}, 
$$

$$
B_{\frac{1}{2}}(n)=\cup\{F\in{\cal{Q}(C)}\,|\,F \cap K_{\frac{1}{2}}(n)\neq\emptyset\} 
$$
and  $F_{0}(n)= p^{-1}(n)\cap J'$. By construction $B_{-\frac{1}{2}}(n)= K_{-\frac{1}{2}}(n)\times [0,1]$
and $B_{\frac{1}{2}}(n)= K_{\frac{1}{2}}(n)\times [0,1]$.\\

Let $B=B_{-\frac{1}{2}}(n)\cup F_{0}(n)\cup B_{\frac{1}{2}}(n)$. Notice that
$B=\mbox{Cl}(p^{-1}(n-1,n+1)\cap J')$, where Cl denotes closure.

\begin{lem}\label{circle}
 The space $B$ is homeomorphic to $\mathbb{S}^{1}\times I$.
\end{lem}

\noindent{\it Proof.} By construction, $B$ is a compact connected submanifold of $J'$. 
Since $J'$ is homeomorphic to $\mathbb{S}^{1}\times\mathbb{R}$, and $J'-B$ has two connected components, then the result follows. $\square$ 

\begin{lem}
 $F_{0}(n)$ has the homotopy type of $S^{1}$.
\end{lem}

\noindent{\it Proof.} Consider the set $B$. By Lemma  \ref{circle}, we know that $B\cong \mathbb{S}^{1}\times I$, where
$K_{-\frac{1}{2}}(n)\times\{0\}\cong \mathbb{S}^1\times\{0\}$ and $K_{\frac{1}{2}}(n)\times\{0\}\cong \mathbb{S}^1\times\{1\}$. 
Hence $\tilde{B}=B/(K_{-\frac{1}{2}}(n)\times\{0\}\cup K_{\frac{1}{2}}(n)\times\{1\})$ is homeomorphic to $\mathbb{S}^2$. By Alexander duality, using reduced homology and cohomology groups,
we have that $\tilde{H}_0(\tilde{B}-F_{0}(n),\mathbb{Z})=\tilde{H}^{1} (F_{0}(n),\mathbb{Z})$. One has that $\tilde{B}-F_{0}(n)$ has two connected components, so
$\tilde{H}_{0} (\tilde{B}-F_{0}(n),\mathbb{Z})\cong \mathbb{Z}$. Since $F_{0}(n)\subset B\cong \mathbb{S}^{1}\times I$, we have that either
$\Pi_1 (F_{0}(n))\cong \{0\}$ or  $\Pi_1 (F_{0}(n))\cong \mathbb{Z}$, but $H^{1}(F_{0}(n),\mathbb{Z})\cong\mathbb{Z}$, hence $\Pi_1 (F_{0}(n))\cong \mathbb{Z}$. 
Therefore, $F_{0}(n)$ has the homotopy type of $S^1$. $\square$\\

Next, we are going to describe the subset $F_{0}(n)$. Notice that the edges of $F_{0}(n)$ are of four types, which we will denote 
by $T_1$, $T_2$, $T_3$ and $T_4$.

\begin{itemize}
\item An edge $l\subset F_{0}(n)$ belongs to $T_1$ if  $l\subset B_{-\frac{1}{2}}(n)$ but $l\not\subset B_{\frac{1}{2}}(n)$.

\item An edge $l\subset F_{0}(n)$ belongs to $T_2$ if $l\subset B_{\frac{1}{2}}(n)$ but $l\not\subset B_{-\frac{1}{2}}(n)$.

\item An edge $l\subset F_{0}(n)$ belongs to $T_3$ if $l\subset B_{-\frac{1}{2}}(n)\cap B_{\frac{1}{2}}(n)$.

\item An edge $l\subset F_{0}(n)$ belongs to $T_4$ if $l\not\subset  B_{\frac{1}{2}}(n)$ and $l\not\subset  B_{-\frac{1}{2}}(n)$.
\end{itemize}
\begin{lem}
The space $B$ retracts strongly to $F_{0}(n)$.
\end{lem}

\noindent{\it Proof.} Since $B=B_{-\frac{1}{2}}(n)\cup F_{0}(n)\cup B_{\frac{1}{2}}(n)$ is homeomorphic to $\mathbb{S}^1\times I$, 
and $B_{-\frac{1}{2}}(n)= K_{-\frac{1}{2}}(n)\times [0,1]$ and $B_{\frac{1}{2}}(n)= K_{\frac{1}{2}}(n)\times [0,1]$, we have that 
$B_{-\frac{1}{2}}(n)= K_{-\frac{1}{2}}(n)\times [0,1]$ retracts strongly to $K_{-\frac{1}{2}}(n)\times \{1\}$ and
$B_{\frac{1}{2}}(n)= K_{\frac{1}{2}}(n)\times [0,1]$ retracts strongly to $K_{\frac{1}{2}}(n)\times \{0\}$. Now
$\partial F_{0}(n)=K_{-\frac{1}{2}}(n)\times\{1\}\cup K_{\frac{1}{2}}(n)\times\{0\}$. Therefore, the result follows. $\square$\\

By the above, we have copies of $K_{-\frac{1}{2}}(n)$ and $K_{\frac{1}{2}}(n)$ contained in $\partial F_{0}(n)$. By abuse of notation we will
denote them in the same way. Notice that $K_{-\frac{1}{2}}(n)$ is the union of edges of types $T_1$ and $T_3$, and
$K_{\frac{1}{2}}(n)$ is the union of edges of types $T_2$ and $T_3$.

\begin{lem}\label{fns}
There exists a finite sequence of cubulated moves that carries the  knot $K_{-\frac{1}{2}}(n)$ into the knot $K_{\frac{1}{2}}(n)$.
\end{lem}

\noindent{\it Proof.} We will show it by cases.\\
\noindent Case 1. Suppose that $K_{-\frac{1}{2}}(n)= K_{\frac{1}{2}}(n)$. The result is obviously true.\\

\noindent Case 2. Suppose that $K_{-\frac{1}{2}}(n)\cap K_{\frac{1}{2}}(n)=\emptyset$. In other words, $K_{-\frac{1}{2}}(n)$ and $K_{\frac{1}{2}}(n)$
 do not have edges of type $T_3$. \\
Remember that $F_{0}(n)$ is a cubulated compact surface whose fundamental group is isomorphic to $\mathbb{Z}$. Thus  $\partial F_{0}(n)$
has two connected component; namely $K_{-\frac{1}{2}}(n)$ and $K_{\frac{1}{2}}(n)$, such that their intersection is empty.
Hence $F_{0}(n)$ is divided by squares (2-faces) belonging to the 2-skeleton of $\cal C$, whose edges are of any of the types $T_1$, $T_2$ and $T_4$.\\

Next, we will carry the knot $K_{-\frac{1}{2}}(n)$ onto the knot $K_{\frac{1}{2}}(n)$ via a finite number of cubulated moves; {\it i.e.} we will
carry the edges of type $T_1$ onto the edges of type $T_2$. Let $F$ be a square contained in $F_{0}(n)$. We can assume, up to an (M1)-move, that if an edge
$l\subset F$ belongs to $T_1$, then $F\cap K_{\frac{1}{2}}(n)=\emptyset$  and
$F\cap K_{-\frac{1}{2}}(n)$  consists of either an edge or two neighboring edges. Analogously, if
 $l\subset F$ belongs to $T_2$, then $F\cap K_{-\frac{1}{2}}(n)=\emptyset$ and
$F\cap K_{\frac{1}{2}}(n)$ consists of either an edge or two neighboring edges.\\

Since $F_{0}(n)$ is compact, then it is the union of a finite number of squares, say $m$. We will enumerate them in the following way. The first square
$F_1$ contains an edge of type $T_1$, the squares $F_n$ and $F_{n+1}$ share an edge $l_n$, and whenever 
it is possible, we choose $F_{n+1}$ in such a way that $l_n$ is parallel to $l_{n-1}$. \\

We will use induction on $m$. Consider $F_1$. We apply an (M2)-move to $F_1$ in such a way that the edges of type $T_4$ are replaced by edges of type  $T_1$.
We consider $F_2$. Observe that $F_1$ and $F_2$ share an edge of type $T_1$. Then we apply again an (M2)-move to replace the edges of type 
$T_4$ by edges of type  $T_1$. We continue inductively. Notice that if $l\subset  F_{0}(n)$ is an edge of type $T_2$, then $l\subset \partial F_{0}(n)$;
so $l$ is not a common edge of two squares $F_i$ and $F_j$ in $F_{0}(n)$; hence if  $l$ is replaced by an edge of type $T_1$, then this replacement is kept
at further steps. Therefore, the result follows.\\

\noindent Case 3. Suppose that $K_{-\frac{1}{2}}(n)\cap K_{\frac{1}{2}}(n)$ consists of a finite number of points.\\
The surface $F_{0}(n)$ consists of connected components $C_i$, $i=1,\ldots,r$ such that each $C_i$ is the union of squares 
$F_{i_1}, \ldots, F_{i_{m_i}}\in{\cal{C}}$ and the intersection $F_i\cap F_j$ is either empty or a point belonging to 
$K_{-\frac{1}{2}}(n)\cap K_{\frac{1}{2}}(n)$. Therefore, we apply the previous argument to each $C_i$.\\

\noindent Case 4. Suppose that the intersection $K_{-\frac{1}{2}}(n)\cap K_{\frac{1}{2}}(n)$ contains an edge of type $T_3$. \\
The surface $F_{0}(n)$ consists of 2-dimensional connected components $C_i$, $i=1,\ldots,r$ and paths $\gamma_{ij}$, where each $C_i$ is a union of
 squares $F_{i_1}, \ldots, F_{i_{m_i}}\in{\cal{C}}$, and $\gamma_{ij}$ is a path (or vertex) joining the component $C_i$ with the component $C_j$.
Observe that if  $\gamma_{ij}$ is a path then it is the union of edges of type $T_3$, hence $\gamma_{ij}\subset K_{-\frac{1}{2}}(n)\cap K_{\frac{1}{2}}(n)$.
Moreover $K_{-\frac{1}{2}}(n)\cap K_{\frac{1}{2}}(n)=\cup \gamma_{ij}$ and  $\partial F_{0}(n)=K_{-\frac{1}{2}}(n)\cup K_{\frac{1}{2}}(n)$.\\

Since $\Pi_1(F_{0}(n))\cong\mathbb{Z}$, then  $C_i$ has the homotopy type either the circle or the disk. 
Suppose that $C_i$ has the homotopy type of the circle,
then $\partial C_i\neq \emptyset$ but $\partial C_i\subset \partial F_{0}(n)$, hence $\partial C_i=\partial F_{0}(n)=K_{-\frac{1}{2}}(n)\cup K_{\frac{1}{2}}(n)$, so
$\partial C_i$ contains an edge $l$ of type $T_3$. However $l$ does not belong to
any square $F$ of $F_{0}(n)$; so $l$ does not belong to $C_i$. This is a contradiction, hence  $C_i$ has the homotopy type of the disk.\\

By the above, $\partial C_i$ is homeomorphic to $\mathbb{S}^1$ and consists of edges of type $T_1$ and $T_2$. Moreover, $\partial C_i$ consists of two arcs
$A_1$ and $A_2$, such that $A_1$ is the union of edges of type $T_1$ ($A_1\subset K_{-\frac{1}{2}}(n)$) and $A_2$ is the union of edges of type $T_2$ ($A_2\subset K_{\frac{1}{2}}(n)$). 
Now we apply the same argument used on case 2, so $A_2$ is replaced by $A_1$. Since we have a finite number of components $C_i$,
the result follows. $\square$\\

\begin{cubemoves_thm} There exists a finite sequence
of cubulated moves that carries $K'_{1}$ into  $K'_{2}$. In other words, $K'_{1}$ is equivalent to  $K'_{2}$ by cubulated moves.\end{cubemoves_thm}

\noindent{\it Proof of Theorem \ref{cubemoves}.}  Recall that there exist integer numbers 
$m_1$ and $m_2$ such that $p^{-1}(t)\cap J' = K'_1$ for all $t\leq m_1$ and
$p^{-1}(t)\cap J' = K'_2$ for all $t\geq m_2$. Consider the integer $m_1+1$. By 
Lemma \ref{fns} there exists a finite number of cubulated moves that
carries the knot $K'_{1}$ into the knot $K_{\frac{1}{2}}(m_1+1)$. We continue inductively, and again by 
Lemma \ref{fns}
there exists a finite number of cubulated moves that
carries the knot $K_{\frac{1}{2}}(m_2-1)$ into the knot $K'_{2}$. Since, we have a finite number of integers on $[m_1,m_2]$, then
there exists a finite sequence of cubulated moves that carries $K'_{1}$ into  $K'_{2}$. $\square$

\section{Discrete knots}

In this section we will use the fact that every knot is isotopic to a cubic knot 
to develop a discrete version of a knot. Since a cubic knot is given by a sequence of edges whose boundaries are 
in the canonical lattice of points with integer coefficients in $\mathbb R^3$,  {\emph i.e.}, the abelian group $\mathbb Z^3$, each knot is determined by a 
cyclic permutation $(a_1,\dots,a_n)$ (with some restrictions), $a_i\in\mathbb Z^3$. 
We will describe a regular diagram of the knot in terms of such cyclic permutations by projecting onto a plane $P$, such that 
it is injective when restricted to the $\mathbb{Z}^3$-lattice and the
image of the $\mathbb{Z}^3$-lattice, $\Lambda_P$, is dense. More precisely, the projection of each knot is determined by a cyclic permutation
$(w_1,\dots,w_n)$ (with some restrictions), $w_i\in \Lambda_P$.\\

Let $\mathcal{C}$ be the standard cubulation of $\mathbb{R}^3$ and let $\cal{S}$ be the corresponding 1-skeleton of $\mathcal{C}$ (scaffolding), which by definition is
determined by the vertices in the lattice $\mathbb{Z}^3$ consisting of point with integer coeffcients, and by edges which are contained in straight lines parallel to the 
coordinate axis and passing through points in the lattice. These straight lines belong to three families of parallel lines: 
The family ${\cal{F}}_1$ of lines parallel to the $x$-axis, the family
${\cal{F}}_2$ of lines parallel to the $y$-axis and the family ${\cal{F}}_3$ of lines parallel to the $z$-axis.

\subsection{Discrete description of cubic knots}

Given a cubic knot $K\subset\mathbb{R}^3$, we have that $K$ is a polygonal simple curve whose vertices lie on the
lattice $\mathbb{Z}^3$ and whose edges lie on straight lines belonging to the families ${\cal{F}}_i$, $i=1,2,3$ of straight lines
parallel to the coordinate axis $x_i$. This implies that
we can write $K$ as a {\it cyclic permutation of points}
$(v_1,v_2,\,\ldots\,,v_n)$ such that $v_i\in \mathbb{Z}^3$, $v_i\neq v_j$ $1\leq i,j\leq n$ and $v_i$ 
is joined to $v_{i+1}$ by a unit edge and $v_{n}$ is likewise joined to $v_{1}$ by a unit edge. 
Observe that the expression  $(v_1,v_2,\,\ldots\,,v_n)$ represents a equivalence class, since
 $(v_1,v_2,\,\ldots\,,v_n)=(v_2,v_3,\,\ldots\,,v_n,v_1)=\cdots=(v_n,v_1,\,\ldots\,,v_{n-1})$.
Conversely, a cyclic permutation of points $(v_1,v_2,\,\ldots\,,v_n)$ satisfying the above conditions 
determines a unique knot $K\subset {\cal{S}}$. Notice that if $K$ is oriented, we can associate to it a unique cyclic permutation
such that the numbering of the $v_i's$ is compatible with this orientation. If $K$ is not oriented there are two cycles
that can be associated to $K$, $(v_1,v_2,\,\ldots\,,v_n)$  
and its inverse $(v_n,v_{n-1},\,\ldots\,,v_1)$ \\ 

Let $e_1,\,e_2,\,e_3$ be the canonical coordinate vectors in $\mathbb{R}^3$. We can also write the cubic knot $K$ as an {\it ``anchored cycle''} permutation of directions
$K=\{v_1,(e_{i_1},\,\ldots\,,e_{i_n})\}$. In other words, 
cubic knots can be coded as  words over
the alphabet $\{e_1^{\pm},e_2^{\pm},e_3^{\pm}\}$ plus the initial vertex.\\

We now can describe the cubulated moves in $\mathbb{R}^3$ using cycle notation.

\begin{itemize}
\item [\bf M1] Subdivision: We subdivide each edge of $K$ by $m$ equal edges. Thus
$ \{v_1,(e_{i_1},\,\ldots\,,e_{i_n})\}$
is equivalent to 
$$
\{v_1,(\frac{1}{m} e_{i_1},\,\ldots,\,\frac{1}{m} e_{i_1},\,\ldots\,,\frac{1}{m} e_{i_n},\,\ldots\,,\frac{1}{m} e_{i_n})\},
$$
where each vector $\frac{1}{m}e_{i_k}$ appears $m$ times (see Figure \ref{FM3}).
\begin{figure}[h] 
 \begin{center}
 \includegraphics[height=3cm]{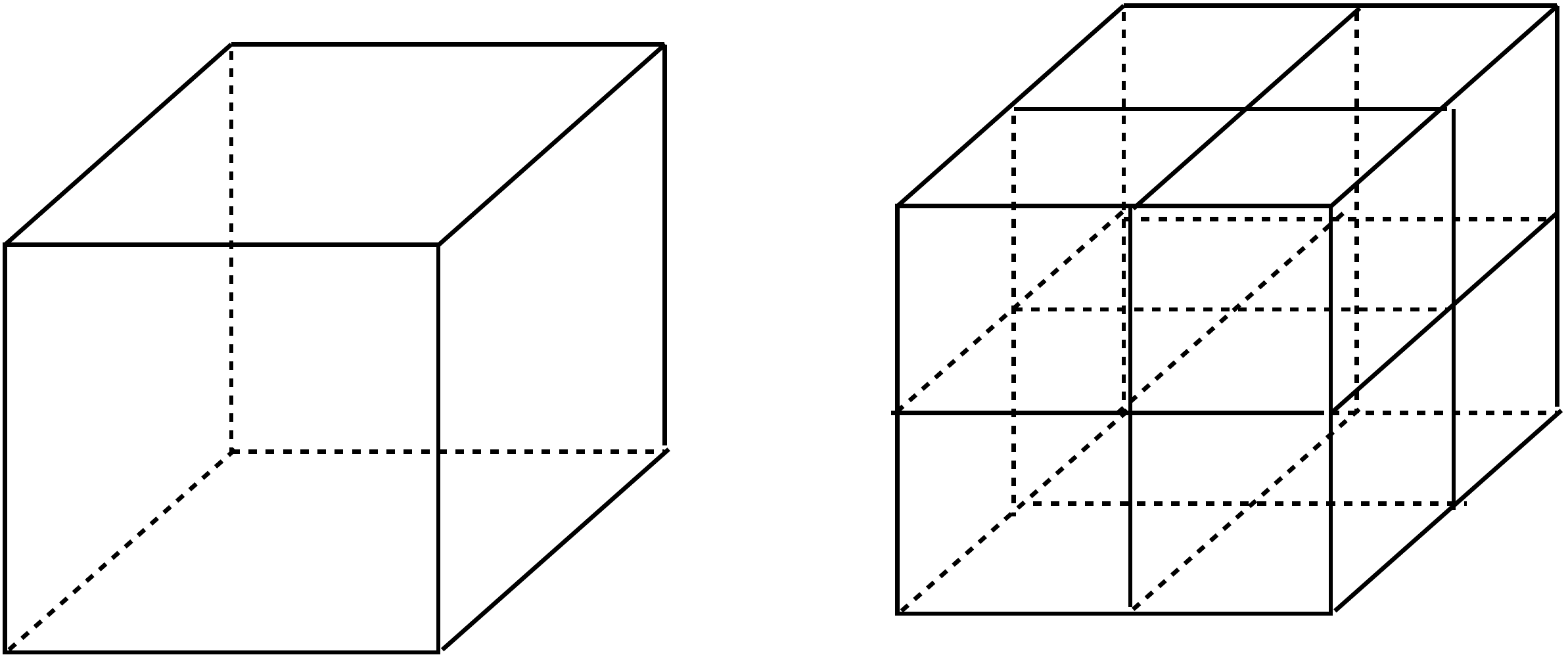}
\end{center}
\caption{\sl Subdivision move by $\frac{1}{2}$.} 
\label{FM3}
\end{figure}

\item[\bf M2] Face boundary move: Let $F$ be a 2-dimensional unit square. Remember that the move M2 consists in replacing the arc $A\subset\partial F$ by 
the arc $A'=\partial F -A$.
If we write $K$ as an anchored cyclic permutation of directions, this move is  expressed as follows:
\begin{enumerate}
\item If the anchored cyclic permutation of directions is given by
$$
\{v_1,(e_{i_1},\,\ldots,\,e_{i_{k-1}},{\bf e_{i_k},e_{i_{k+1}},-e_{i_k}},e_{i_{k+3}},\,\ldots,\,,e_{i_n})\}
$$ 
and $\{v_1,(e_{i_1},\,\ldots\,,e_{i_{k-1}},{\bf e_{i_{k+1}}},e_{i_{k+3}},\,\ldots\,,e_{i_n})\}$\\
are equivalent (see Figure \ref{FM1}). 

In addition, we have: 
$$
\{{\bf v_1}, ({\bf e_{i_1}},e_{i_2}\,\ldots,\,e_{i_{n-2}},{\bf -e_{i_1},e_{i_{n}}})\}
$$ 
and $\{{\bf v_1+e_{i_1}}, (e_{2},\,\ldots,\,e_{i_{n-2}},{\bf e_{i_{n}}})\}$ are equivalent, and 
$$
\{{\bf v_1}, ({\bf e_{i_1}},e_{i_2},e_{i_3},\,\ldots,\,e_{i_{n-1}},{\bf -e_{i_2}})\} 
$$ 
and $\{{\bf v_1+e_{i_2}}, ({\bf e_{1}},e_{i_3},\,\ldots,\,e_{i_{n-2}},{\bf e_{i_{n}}})\}$ are equivalent.

\begin{figure}[h] 
 \begin{center}
 \includegraphics[height=3cm]{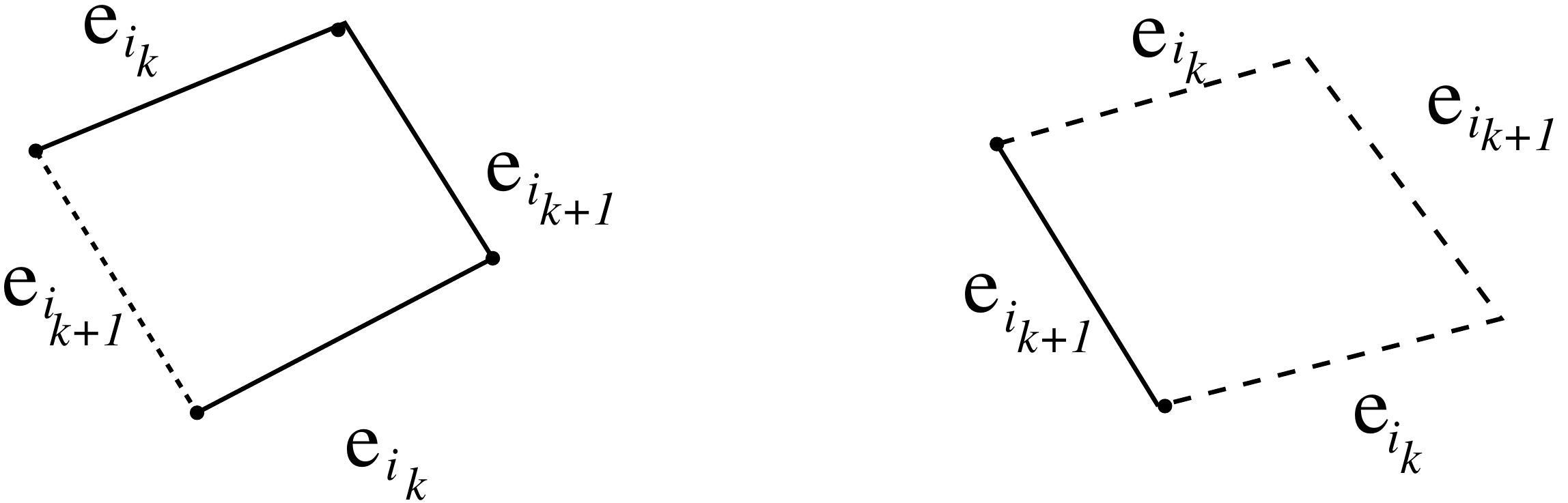}
\end{center}
\caption{\sl Face boundary move.} 
\label{FM1}
\end{figure}

\item If the anchored cyclic permutation of directions is given by
$$
\{v_1,(e_{i_1},\,\ldots\,,e_{i_{k-1}},{\bf e_{i_k},e_{i_{k+1}}},e_{i_{k+2}},\,\ldots\,,e_{i_n})\}
$$ 
and $\{v_1,(e_{i_1},\,\ldots\,,e_{i_{k-1}},{\bf e_{i_{k+1}},e_{i_k}},e_{i_{k+2}},\,\ldots\,,e_{i_n})\}$ are equivalent
(see Figure \ref{FM2}). 

In addition:
$$
\{{\bf v_1, (e_{i_1}},e_{i_2}\,\ldots\,,e_{i_{n-1}},{\bf e_{i_n}})\}
$$ 
and $\{{\bf v_1+e_{i_1}-e_{i_n}}, (e_{i_n},e_{i_2}\,\ldots\,,e_{i_{n-1}},{\bf e_{i_1}})\}$ are equivalent.

\begin{figure}[h] 
 \begin{center}
 \includegraphics[height=3cm]{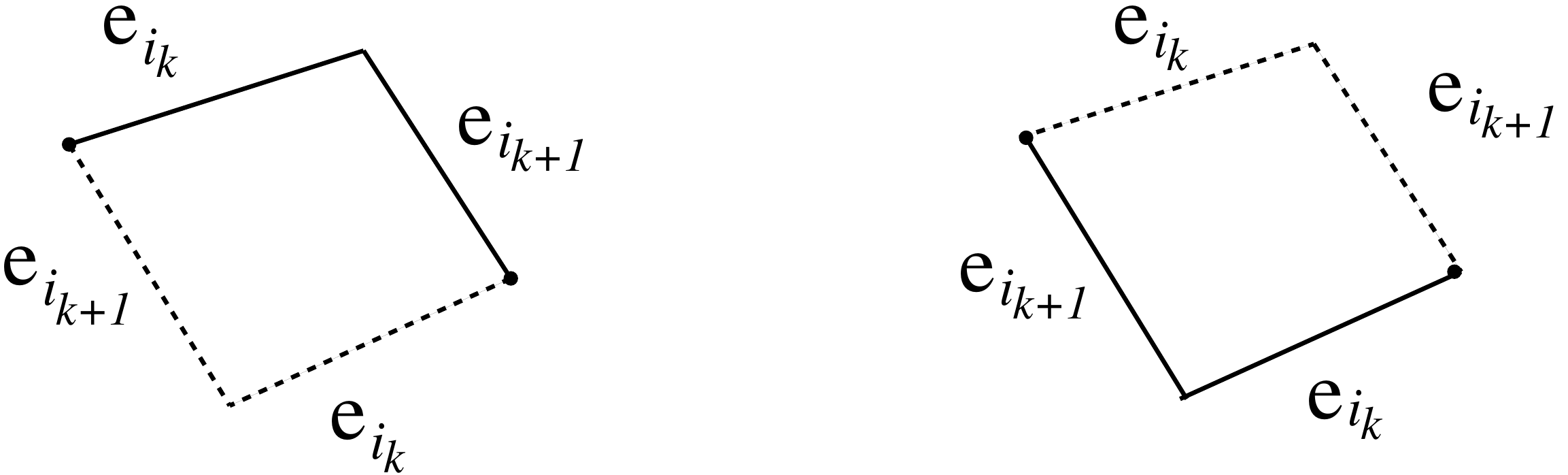}
\end{center}
\caption{\sl Face boundary move.} 
\label{FM2}
\end{figure}
\end{enumerate}
\end{itemize}

\begin{rem}
Note that the face boundary or M2-move in case 2, is only allowed if the vertex  that replaces $v_{k+1}$ or $v_{i_1}$, 
{\it i.e.}, $v_{i_{k+1}}+e_{i_{k+1}}$ or $v_1+e_{i_1}-e_{i_n}$ is not already a vertex of the knot. However, one can always apply the subdivision or M1-move in that case.
\end{rem}

\begin{definition}
 Given two anchored cyclic permutation of directions\\
 $\{v_1,(e_{i_1},\,\ldots\,,e_{i_n})\}$, $\{b_1,(e_{j_1},\,\ldots\,,e_{j_r})\}$, we say that they are
{\emph{equivalent}} if there exists a  finite sequence of cubulated moves 
transforming the anchored cyclic
permutation of directions $\{h(v_1),(e_{i_1},\,\ldots\,,e_{i_n})\}$ into $\{b_1,(e_{j_1},\,\ldots\,,e_{j_r})\}$, where $h$ is the translation map
given by $h(x)= x+b_1-v_1$.
\end{definition}

By Theorem \ref{main},  we know that given two cubic knots $K_1$ and $K_2$  in $\mathbb{R}^{3}$, then
they are isotopic if and only if there exists a finite sequence of cubulated moves that carries $K_1$ into $K_2$. If we express both knots
as  anchored cyclic permutations of directions $K_1=(v_1,e_{i_1},\,\ldots\,,e_{i_n})$, $K_2=(v_2,e_{j_1},\,\ldots\,,e_{j_r})$, we have the following

\begin{coro}\label{eqperm}
Two cubic knots $K_1$ and $K_2$  in $\mathbb{R}^{3}$ are isotopic if and only if their anchored cyclic permutations of directions 
$K_1=\{v_1,(e_{i_1},\,\ldots\,,e_{i_n})\}$ and $K_2=\{b_1,e_{j_1},\,\ldots\,,e_{j_r})\}$
are equivalent.
\end{coro}

\subsection{Discrete description of a suitable cubic knot projection}

Many invariants of a knot are computed using suitable projections. We will start proving that for any
 oriented cubic knot $K=(v_1,v_2,\,\ldots\,,v_n)$, there exists  a generic projection $p$ over a suitable fixed plane $P$. Next, we will represent $\widehat{K}= p(K)$
 as a cyclic permutation of points $(w_1,w_2,$ $\,\ldots\,,w_n)$.

\subsubsection{A canonical generic projection}

\begin{lem}\label{inyectiva}
There exists a hyperplane $P$ such that the orthogonal projection $p:\mathbb{R}^3\rightarrow P$ is injective in the vertices of  $\cal{S}$ (the lattice $\mathbb{Z}^3$) and is also 
injective in each of the families ${\cal{F}}_i$, $i=1,2,3$. 
\end{lem}

\noindent{\it Proof.} Let $N=(1,\pi,\pi^2)$, where $\pi$ is the well-known transcendental number.
 Let $P$ be the  plane  through the origin in $\mathbb{R}^3$ orthogonal to $N$.
 Let us show that $p:\mathbb{R}^3\rightarrow P$ is injective when restricted to the lattice $\mathbb{Z}^3$.
Let $(x_0,y_0,z_0)$ be a point in $P$ and $L=\{t(1,\pi,\pi^2)+(x_0,y_0,z_0):t \in \mathbb{R}\}$ the line which is orthogonal to $P$ and passes through $(x_0,y_0,z_0)$.
Let us suppose that $L$ contains two points with integer coefficients, $P_1=(N,M,K)=(x_0 + t_1, y_0 +\pi(t_1), z_0 + \pi^2(t_1))$ and
 $P_2=(N',M',K')=(x_0 + t_2, y_0 +\pi(t_2), z_0 + \pi^2(t_2))$. If we consider the point $P_1-P_2$, which belongs to $L$, we must have that
$P_1=P_2$ since otherwise $\pi$ and $\pi^2$ would be algebraic numbers which is a contradiction since $\pi$ is transcendental.
Therefore $p$ restricted to the lattice $\mathbb{Z}^3$ is
injective. 
Now we will prove that the projection $p$ restricted to the family ${\cal{F}}_i$, $i=1,2,3$ is injective. 
First we will consider the case $i=1$. Let $l$ and $l'$ be two lines belonging to ${\cal{F}}_1$. 
Suppose that $p(l)=p(l')$. Let $Q$ be the plane in $\mathbb{R}^3$ containing both $l$ and $l'$. Since $p$ restricted to the coordinate planes is injective,
we can assume that $Q$ is not parallel to a coordinate plane. Therefore, a normal vector to $Q$ is of the form $(0,a,b)$ where both $a$ and $b$ are nonzero. 
Also  $p(Q)$ is the line  $p(l)=p(l')$. 
This means that $Q$ contains a translate of the vector $(1,\pi,\pi^2)$ normal to $P$. That is, $(1,\pi,\pi^2)+x$ belongs to $Q$, (here $x$ is a point in $Q$).
Hence the interior product of $(1,\pi,\pi^2)$ and $(0,a,b)$ is zero, which is impossible because it would mean that $\pi$ is a rational multiple of $\pi^2$. 
By a similar argument $p$  restricted to each of the families ${\cal{F}}_2$ and ${\cal{F}}_3$ is  also injective. $\square$\\

\begin{lem}
For $p:\mathbb{R}^3\rightarrow P$ defined as above,  $p(\mathbb{Z}^3) := \Lambda_P$ is an additive dense subgroup of $P$.
Furthermore, the images $p({\cal{F}}_i)$ of each of the families ${\cal{F}}_i$, $i=1,2,3$ are dense and invariant under the group of translations of the dense group
$\Lambda_P$. In addition, for  $i=1,2,3$, $p({\cal{F}}_i)$ is a one dimensional dense subgroup of the additive Lie group $({\mathbb R}^3,+)$.
\end{lem}

\noindent{\it Proof.} Let $e_1,\,e_2,\,e_3$ be the canonical coordinate vectors in $\mathbb{R}^3$ and $f_i=p(e_i),\,\,i=1,2,3$ be their projections in $P$. Since $p$ restricted
to the coordinate lines is injective, we have that $f_1,\,f_2$ and $f_3$ are pairwise linearly independent. Consider the lattice $\mu=af_1+bf_2$,  $a,b\in\mathbb{Z}$.
It is enough to prove that the points $kf_3$, $k\in\mathbb{Z}$ are dense in the torus $P/\mu$. Since $f_3=\alpha f_1+\beta f_2$,  $\alpha,\,\beta\in\mathbb{R}$, we will prove that 
both $\alpha$ and $\beta$ are irrational numbers. Suppose that  $f_3=\alpha f_1+\frac{n}{m} f_2$. Then $p(mf_1-nf_2-m\alpha f_3)=0$, {\it i.e.},
$(m,-n,-m\alpha)\in\ker(p)$. As $\ker(p)$ is the line orthogonal to $P$, $L_N=\{t(1,\pi,\pi^2)\,:\,t \in \mathbb{R}\}$, it follows that 
$(m,-n,-m\alpha)= t(1,\pi,\pi^2)$  for some $t \in \mathbb{R}$. This implies that $\pi$ is rational, which is a contradiction.

\begin{rem} The set $\Lambda_P$
is a dense subgroup of rank 3 of $P$. In fact, if $e_1,\,e_2,\,e_3$ denote the canonical coordinate vectors,
then $\mathbb{Z}^3=\{m_1e_1+m_2e_2+m_3e_3\,|\,m_1,m_2,m_3\in\mathbb{Z}\}$ is an additive subgruop of $\mathbb{R}^3$ and since the projection
$p$ is a group homomorphism, we have that  $\Lambda_P=\{m_1f_1+m_2f_2+m_3f_3\,|\,m_i\in\mathbb{Z};\,\,f_j=p(e_j)\,\,j=1,2,3\} $ is a free
abelian group with three generators.\\
\end{rem}

\noindent{\large \bf Acknowledgments}\\

\noindent The authors thank Professor Margareta Boege for her helpful suggetions and comments. We would also like to thank Professor Scott Baldridge for pointing out the relationship between
our work and his paper. Gabriela Hinojosa was partially supported by CONACyT (M\'exico), 
CB-2009-129939 and Alberto Verjosvky was partially supported by
CONACyT (M\'exico), CB-2009-129280 and PAPIIT (Universidad Nacional Aut\'onoma de M\'exico) IN100811.

G. Hinojosa. {\tt Facultad de Ciencias}. Universidad Aut\'onoma del Estado de Morelos. Av. Universidad 1001, Col. Chamilpa.
Cuernavaca, Morelos, M\'exico, 62209.

\noindent {\it E-mail address:} gabriela@uaem.mx
\vskip .3cm
A. Verjovsky. {\tt Instituto de Matem\'aticas, Unidad Cuernavaca}. Universidad Nacional Au\-t\'o\-no\-ma de M\'exico.
Av. Universidad s/n, Col. Lomas de Chamilpa. Cuernavaca, Morelos, M\'exico, 62209.

\noindent {\it E-mail address:} alberto@matcuer.unam.mx
\vskip .3cm
C. Verjovsky Marcotte. {\tt Mathematics Department}.  St. Edwards University. 3001 South Congress Avenue, Austin, TX 78704, USA.

\noindent {\it E-mail address:} cynthiavm@gmail.com

\end{document}